\input amstex
\input epsf
\documentstyle{amsppt}
\nologo

\def\proof{\demo{\smc Proof}}
\def\endproof{\hfill$\square$\enddemo}
\overfullrule=0pt

\topmatter
\title Solving the Selesnick-Burrus Filter Design Equations
Using Computational Algebra and Algebraic Geometry\endtitle
\leftheadtext{Solving the Selesnick-Burrus Filter Design Equations}
\rightheadtext{Solving the Selesnick-Burrus Filter Design Equations}
\author John B. Little\endauthor
\affil Department of Mathematics and Computer Science\\
       College of the Holy Cross \endaffil
\address Department of Mathematics and Computer Science
         College of the Holy Cross
         Worcester, MA 01610\endaddress
\date September 19, 2002\enddate
\email   little\@mathcs.holycross.edu \endemail
\abstract 
In a recent paper, I. Selesnick and C.S. Burrus developed a design method
for maximally flat FIR low-pass digital filters with 
reduced group delay.  Their approach leads to a system of 
polynomial equations depending on three integer design 
parameters $K,L,M$.  In certain cases (their ``Region I''), 
Selesnick and Burrus were able to derive solutions using only 
linear algebra; for the remaining cases (``Region II''), they 
proposed using Gr\"obner bases.  This paper introduces a different method,
based on multipolynomial resultants, for analyzing and solving the 
Selesnick-Burrus design equations.  The results of calculations
are presented, and some patterns concerning the number of solutions 
as a function of the design parameters are proved.
\endabstract
\endtopmatter

\document

\heading \S 1. Introduction \endheading

In this paper we will present an application of techniques from
computational commutative algebra and algebraic geometry to a
problem in signal processing.  We will see that recent
developments in the theory of multipolynomial resultants give
a powerful method for solving an interesting family of problems 
in digital filter design.

We begin by recalling some basic concepts about digital filters.  
(A good general reference for this material is \cite{PM}.)
A digital signal is a quantized function of a 
discrete variable, (e.g. time). If we ignore quantization
effects, therefore, a signal can be represented mathematically 
by a sequence of complex numbers $x[n]$ indexed by $n\in {\Bbb Z}$.
For many purposes, an appropriate class of signals is the sequence
space $\ell_2$, since the finiteness of the $\ell_2$ norm corresponds
to a finite energy condition on signals.
Signal processing operations can be 
described mathematically by means of operators 
$\Gamma : \ell_2 \to \ell_2.$
In the signal processing context, these are called {\it digital filters}.
Here, we only consider filters that are linear and 
shift-invariant:  If $k$ is fixed and $y[n] = x[n + k]$ for all $n$, 
then $\Gamma(y)[n] = \Gamma(x)[n+k]$.  

A linear, shift-invariant
filter is characterized completely by its {\it transfer function} $H(z)$, 
the $z$-transform of its impulse response  (see \S 1 below).
Design methods for filters to perform specified operations on 
signals can often be formulated as finding solutions of systems
of polynomial equations on the coefficients in transfer functions
$H(z)$ of some specified form.  For this 
reason, techniques from computational commutative algebra have begun
to find uses in this area.  

In this article we will focus on one particular filter design
method introduced by Selesnick and Burrus in \cite{SB}.  Their idea
was to specify $H(z)$ for a low-pass, finite impulse response (FIR) 
filter (see \S 1) by imposing three types of conditions:
\roster
\item  A given number $M$ of flatness conditions at $\omega = 0$
on the {\it square magnitude response} 
$$F(\omega) = |H(e^{i\omega})|^2$$
(that is, the vanishing of the derivatives of all orders up to $2M$ 
of $F(\omega)$ at $\omega = 0$ -- note that $F$ is an even function
so the derivatives of odd orders at $\omega = 0$ are zero automatically),
\item  A second number $L$ of flatness conditions at $\omega = 0$
on the {\it group delay} 
$$G(\omega) = {d\over d\omega} \arg H(e^{i\omega}) $$
(that is, the vanishing of the derivatives of all orders up to $2L$ 
of $G(\omega)$ at $\omega = 0$ -- note that $G$ is also an even 
function of $\omega$),
and
\item A third number $K$ of zeroes of $H(e^{i\omega})$ at $\omega = \pi$.
\endroster
The parameters $K,L,M$ can be specified independently and 
this approach can be seen as a generalization of earlier work on 
maximally flat filters by Hermann, Baher, and others in certain special cases.  
Each of these types conditions leads to  polynomial 
equations of degree $\le 2$ 
on the coefficients $h[n]$ in $H(z) = \sum_{n=0}^{N-1} h[n]z^{-n}$, and
solutions exist provided $N - 1 \ge K + L + M$.
The equations have a particularly simple form if the filter moments
$$m_k = \sum_{n = 0}^{N-1} n^k h[n]$$
are used as the variables.  Following Selesnick and Burrus,
we express everything in terms of the $m_k$.

Selesnick and Burrus establish a subdivision of these 
problems into two classes.  The {\it easier} cases (Region I) occur for 
$L$ relatively large compared to $M$:  
$$\left\lfloor {M-1\over 2}\right\rfloor \le L \le M.$$
In these cases, Selesnick and Burrus give an analytic solution procedure
depending only on linear algebra.  
The more {\it difficult} cases (Region II) occur when $L$ is relatively small
compared to $M$:
$$0 \le L \le \left\lfloor {M-1\over 2}\right\rfloor - 1$$
In Region II, Selesnick and Burrus used {\it lex} Gr\"obner basis computations
to solve the resulting filter design equations in a few cases.  However,
the complexity of this approach severely limited the range of cases they
were able to handle.

Some remaining problems left unsolved by Selesnick and Burrus's work are 

\bigskip
\roster
\item to develop an efficient method to solve the filter design 
equations in the Region II cases, and
\item to understand the structure of the 
solutions of the equations for Region II in more detail -- in particular
to determine for given $K,L,M$, how many solutions there are, 
how many are real, how many yield monotone decreasing square
magnitude response $|H(e^{i\omega})|^2$, and so forth.
\endroster

\bigskip
While we cannot claim a complete solution to these problems, 
in this article we first introduce a different solution strategy for the 
Selesnick-Burrus equations in the 
Region II cases which has allowed us to compute solutions
in cases with much larger values of $L,M$ than those reported in \cite{SB}.  
Our approach is based on a careful study of the form of the equations, 
combined with an application of  multipolynomial resultants 
to eliminate variables and obtain a univariate polynomial in the 1st filter 
moment $m_1$.  This strategy is laid out in more detail in (3.10) below.
(For general background on multipolynomial
resultants, see \cite{CLO} Chapters 3 and 7, \cite{EM},
\cite{S} and \cite{CE} for more details on the sparse version, 
and \cite{KSY} for Dixon resultants.  \cite{M} contains
a number of practical recipes for applying these ideas to 
solve systems of equations.)

Second, we attempt to explain some of the intriguing patterns
we have noticed in the solutions, in particular in the number
of distinct complex solutions of the Selesnick-Burrus equations
along the ``diagonals'' $M = 2L + q$ for various values of $q$.
For a given $q$ and $L$ sufficiently large these systems have
a similar shape, and for the first few values of $q$ giving
cases in Region II, we have been able to analyze the form of 
the resultant and determine the degree of the univariate polynomial
in $m_1$ obtained by elimination in all cases.

The organization of the paper is as follows.  \S 2 contains some
additional concepts and terminology on digital filters, 
a presentation of the exact form of the Selesnick-Burrus
equations from \cite{SB}, and a small example (the case $K = 1, L = 1, M = 5$),
which illustrates some key features of these problems.
In \S 3, we lay out a successful solution strategy for the 
Region II problems based on resultants.  The first step consists
of two reductions that permit the direct elimination
of variables in the full Selesnick-Burrus system of $K+L+M$
equations in $K+L+M$ unknowns to yield a much more
manageable system of $M - L - 1$ equations in $M - L - 1$ variables
that we call the reduced Selesnick-Burrus system.  The general
strategy is presented, followed by some experimental results.

First, we present an outline of a calculation determining 
the real solutions of the Selesnick-Burrus system with $K = 2, L = 2, M = 10$,
and the square magnitude response curves of the corresponding filters.,
For this calculation we use a method based on the Dixon resultant, 
combined with numerical rootfinding.  All the calculations were
carried out in the Maple 8 computer algebra system.

Second, we give a table showing the number
of distinct solutions of the Selesnick-Burrus systems for 
most of the cases with $M \le 14$ in Region II (see Figure 2 below).  
A number of the entries in this table were computed by Robert Lewis of 
Fordham University using his {\it Fermat} system and code for Dixon
resultants.  The resultant
strategy would allow the computation of many additional cases with 
$M \ge 15$ as well.  By way of comparison, we note that Selesnick
and Burrus were only able to handle cases with $M \le 7$ in 
their paper.  

In the remainder of the paper we study some of the patterns 
that are apparent in Figure 2.  
\S 4 is devoted to 
a study of the properties of the coefficient matrices 
of the linear parts of the reduced Selesnick-Burrus systems, 
matrices whose coefficients are polynomials in the variable
$t = m_1$.  By some fairly intricate algebraic maneuvering, we
are able to express these matrices in a very useful form
using some notions from the calculus of finite differences.
In particular, the entries can be expressed in terms of  
polynomials of the form $D^j_K(i-t)^\ell|_{i=0}$, 
where $D^j_k$ are certain finite difference operators.  
This allows us to determine the Smith normal form of these
matrices, hence to completely understand the dependence of the
ranks of various submatrices on $t$. 

The cases with $M = 2L + 3$ are studied intensively 
in \S 5, and the following main theorem is established
(compare with the data in the table in Figure 2 below).
\bigskip
\proclaim{(5.1) Theorem} In the cases $M = 2L + 3$, $L\ge 0$, 
(the ``corners'' in Region II boundary), for all $K \ge 1$, the univariate
polynomial in $t$ in the elimination ideal of the Selesnick-Burrus
equations obtained via Strategy (3.10) has degree $8L+8$.
\endproclaim

In \S 6, we undertake a similar study of the 
cases with $M = 2L + 4$ and establish our second main theorem.
\bigskip
\proclaim{(6.1) Theorem} In the cases $M = 2L + 4$, $L\ge 1$, 
for all $K \ge 1$, the univariate polynomial in $t$ in the elimination 
ideal of the Selesnick-Burrus equations obtained via Strategy (3.10) 
has degree $12L + 14$.
\endproclaim

The proofs of Theorems (5.1) and (6.1) show in essence how
to construct the appropriate resultant matrices, so they 
give a general, extremely efficient, way to solve all 
cases with $M = 2L+3, 2L+4$.  Similar results are 
possible in principle for the lower diagonals $M = 2L + q$, 
$q \ge 5$
as well.  But we will not attempt to prove formulas for the
number of solutions in those cases here
because the resultants necessary to handle them become
progressively more complicated to analyze.

In a companion article, \cite{LL}, we will discuss the properties of the
Selesnick-Burrus filters from Region II in more detail.  

The author would like to thank Ivan Selesnick for several valuable
conversations, and Robert Lewis for permission to present his computational 
results here.

\heading \S 2. Preliminaries on Filter Design and the Selesnick-Burrus Equations \endheading

\bigskip

Let $\delta$ be the signal
$$\cdots, 0,0,1,0,0, \cdots $$
($1$ at $n = 0$).  $\delta$ is called the {\it unit impulse} at $n = 0$.
Let $\Gamma$ be a linear, shift-invariant filter as in \S 1.  The 
output $\Gamma(\delta)$ from the filter on input $\delta$ is called the 
{\it impulse response} of $\Gamma$.  A beautiful consequence
of the linearity and shift invariance hypotheses is that 
the impulse respose of a filter determines
the output on any other input signal.  For, we can write
$$x[n] = \sum_{k=-\infty}^\infty x[k]\delta[n-k].$$
If $h[n]$ are the coefficients of the impulse response and
and $y = \Gamma(x)$ is the output, then by linearity and 
shift-invariance, 
$$\eqalign{
y[n] &=  \sum_{k=-\infty}^\infty x[k] \Gamma(\delta)[n-k]\cr
     &= \sum_{k=-\infty}^\infty x[k] h[n-k]\cr}\leqno(2.1)$$
In other words, {\it the output is the (discrete) convolution
of the input and the impulse response}.

It is standard in signal processing to    
package the signals $x[n],y[n],h[n]$ by their ``$z$-transforms''
$X,Y,Z$.  For instance, the definition of the $z$-transform of the
signal $x[n]$ is
$$X(z) = \sum_{n=-\infty}^\infty x[n] z^{-n}.$$

The $z$-transform of the impulse response, $H(z)$, 
is called the {\it transfer function} of the filter.
In our cases, $h[n]$ will be nonzero for only finitely many $n$.
Such filters are called {\it finite impulse response}, or FIR
filters.  For an FIR filter, the transfer function is a rational
function, hence has a well-defined value at all $z$ in the complex
plane, except for a pole at $z = 0$.  

Note that the coefficient of $z^{-n}$ in the product
$H(z)X(z)$ is the discrete convolution from (2.1)
$$\sum_{k=-\infty}^\infty h[k]x[n-k],$$
which is the same as $y[n]$.  In other words, {\it the
$z$-transform of the output is the product of the transfer function
and the $z$-transform of the input:}  $Y(z) = H(z)X(z)$.

Note that the restriction of $H(z)$ to the unit circle in 
the complex plane, 
$$H(e^{i\omega}) = \sum_{n=-\infty}^\infty h[n] e^{-i n\omega},$$
is the (discrete-time) Fourier transform 
of $h$, so $H(z)$ also determines the 
{\it frequency response} characteristics of the filter on 
input signals.

Filter design problems, such as the one studied in \cite{SB}, 
ask for constructions of filters adapted to perform some specified operation
on input signals.  An important approach is to obtain the desired
behavior by designing the form of the transfer function $H(z)$.
For instance, we might seek to construct:

\bigskip
\roster
\item ``Low-pass'' filters in order to remove high-frequency components
of signals.  These typically {\it smooth out} or {\it blur} signals
and can be used to remove high-frequency ``noise''.
\item ``High-pass'' filters to remove low-frequency components of
input signals.  These typically pick out fine details, or rapid changes
in the input and can be used to detect features.
\endroster

\bigskip
The paper of Selesnick and Burrus proposes a way to design
maximally flat low-pass FIR filters with reduced group delay.
These filters are specified by three positive integer parameters denoted
$K,L,M$.  For an FIR low-pass filter with transfer function
$$H(z) = \sum_{n=0}^{N-1} h[n]z^{-n},$$
let $F(\omega)$ be the square magnitude response and 
$G(\omega)$ be the group delay response as in \S 1.
Selesnick and Burrus show that if $K,L,M \in {\Bbb N}$, and
$K + L + M + 1 = N$, $L \le M$, then the filter coefficients
$h[n]$ can be determined to make:

$$\eqalign{F^{(2i)}(0) &= 0, \quad i = 1, \ldots, M,\cr
           G^{(2j)}(0) &= 0, \quad j = 1, \ldots, L,\cr
           (1+z^{-1})^K &\ |\  H(z). \cr}\leqno(2.2)$$

The meaning of the first condition is that $F(\omega)$
is {\it flat to order $2M$ at $\omega = 0$}.  Similarly the
second equation says $G(\omega)$ is flat to order $2L$ at $\omega = 0$.
The final equation can also be interpreted as a flatness condition,
since it implies that $|H(\omega)|^2$ has a zero of order $2K$ at
the normalized frequency $\omega = \pi$, which corresponds to 
$z = -1$ under $z = e^{i\omega}$.

It is easy to see that the Selesnick-Burrus conditions (2.2) 
can be expressed as polynomial equations in the filter coefficients.
However, the form of these equations becomes significantly simpler if 
they are expressed in terms in terms of the filter {\it moments},
$$m_k = \sum_{n=0}^{N-1} n^k h[n].\leqno(2.3)$$
The explicit form of the equations is derived in \cite{SB} as follows:
\bigskip
\item{1.}  The flatness conditions on $F$ at $\omega = 0$
are quadratic in the $m_i$: 
$$0 = {2i\choose i} m_i^2 + 2\sum_{\ell=0}^{i-1} {2i\choose \ell}(-1)^{i+\ell} m_\ell m_{2i-\ell},
\quad i = 1,\ldots, M.\leqno(2.4a)$$
\item{2.} The flatness conditions on $G$ at $\omega = 0$
are also quadratic in the $m_i$: 
$$0 = \sum_{\ell=0}^j \left(1-{2\ell\over 2j+1}\right) {2j+1\choose \ell} (-1)^\ell m_\ell m_{2j+1-\ell},
\quad j = 1,\ldots, L.\leqno(2.4b)$$
(These are derived from $G^{(2j)}(0) = 0$, using the conditions
$F^{(2i)}(0) = 0$, $i = 1, \ldots, M$.)
\item{3.} Finally, the zero of order $K$ at $z = -1$ is equivalent
to saying that the remainder of $H(z)$ on division by $(1+z^{-1})^K$ is
zero.  This yields $K$ linear equations on $m_i$.

\bigskip
At first glance this looks like an underdetermined system with
$2M+1$ variables $m_i$, $i = 0,\ldots, 2M$,
and $K + L + M = N - 1$ equations.  However, the moments $m_k$, $k \ge N$ are not 
independent variables.  They can all be expressed in terms of 
$m_0, \ldots, m_{N-1}$ by solving systems of linear equations.  
We will normalize our filters by requiring that $m_0 = 1$.
This accomplishes a first reduction to a system of $N-1$ equations in $N-1$
variables.  We expect only finitely many solutions and the real solutions
are of the greatest interest.

\bigskip
\noindent
{\bf (2.5) Example.} We study the Selesnick-Burrus equations
in the relatively simple case $L = 1, M = 5, K = 1$.  There are 6 
quadratic equations, from setting
$$\eqalign{
&2 m_1^2 - 2 m_2,\cr
&6 m_2^2 + 2 m_4 - 8 m_1 m_3,\cr
&20 m_3^2 - 2 m_6 + 12 m_1 m_5 - 30 m_2 m_4,\cr
&70 m_4^2 + 2m_8 - 16 m_1 m_7 + 56 m_2 m_6 - 112 m_3 m_5,  \cr
&252 m_5^2 - 2m_{10} + 20 m_1 m_9 - 90m_2 m_8 + 240 m_3 m_7 - 420 m_4 m_6\cr
& -m_3 + m_1 m_2\cr}$$
equal to zero, and similarly 4 additional linear equations:
$$\eqalign{
& - 315 + 14496m_{1} + 23912m_{3} - 9310m_{4} + 8m_{7} - 196m_{6} + 1904m_{5} - 30184m_{2},  \cr
&2m_{8} - 728m_{6} + 9408m_{5} - 51632m_{4} + 141120m_{3} - 185152m_{2} + 91392m_{1} - 2205, \cr
&4m_{9} - 17052m_{6} + 247380m_{5} - 1445010m_{4} + 4105160m_{3} - 5529048m_{2}  \cr
&\quad + 2784096m_{1} - 72765, \cr
&m_{10} - 43407m_{6} + 670320m_{5} - 4070200m_{4}+ 11869200m_{3} - 16288944m_{2} \cr
&\quad + 8326080m_{1} - 231525.\cr
}$$

In this small example, we can apply a ``brute force'' method to 
derive a solution.  This is also essentially the method used by 
Selesnick and Burrus to handle the more difficult problems in their Region II.
 The {\it lex} Gr\"obner basis for the whole system 
with $m_{10} > m_9 > \cdots > m_1$ is in generic ``Shape Lemma'' (\cite{CLO},
Chapter 2, \S 4) form.
The last element is a univariate polynomial of degree 16 in $m_1$.
Using numerical root-finding, we find 6 approximate real roots:
$m_1 \doteq .04470426799$, $1.233505559$, $2.558981682$, $4.441018318$, $5.766494441$, 
and $6.955295732$.  Then the other moments $m_j$ and filter coefficients
$h[i]$ can be determined by backsolving in the Gr\"obner basis and using 
the equations (2.3).  

We can see a general feature of the Selesnick-Burrus
equations here.  Note that the 6 real roots form three pairs of the form $r,7-r$. 
In fact, for all $K,L,M$, the mapping 
$$m_1 \mapsto (L+M+K) - m_1\leqno(2.6)$$ 
gives the effect of {\it time reversal}
(that is, taking the original transfer function 
$H(z) = \sum_{n=0}^{N-1} h[n]z^{-n}$ to the reversed  
$\tilde{H}(z) = \sum_{n=0}^{N-1} h[N-1-n]z^{-n}$).  It is not difficult
to see that the whole Selesnick-Burrus system -- (2.4a), (2.4b), and the 
linear equations expressing the higher moments in terms of the lower ones --
is invariant under time reversal.
Up to time reversal, there are 3 distinct real filters satisfying the 
Selesnick-Burrus conditions in this case.  The plot in Figure 1. shows the 
square magnitude response curves for the three filters.  Note
that two are apparently monotone decreasing, while one has a pronounced
``ripple'' in the ``passband''.   The filters with monotone square magnitude 
responses would be much more useful for actual low-pass filtering 
applications.
\bigskip

\epsfxsize 1.5truein
\centerline{\epsfbox{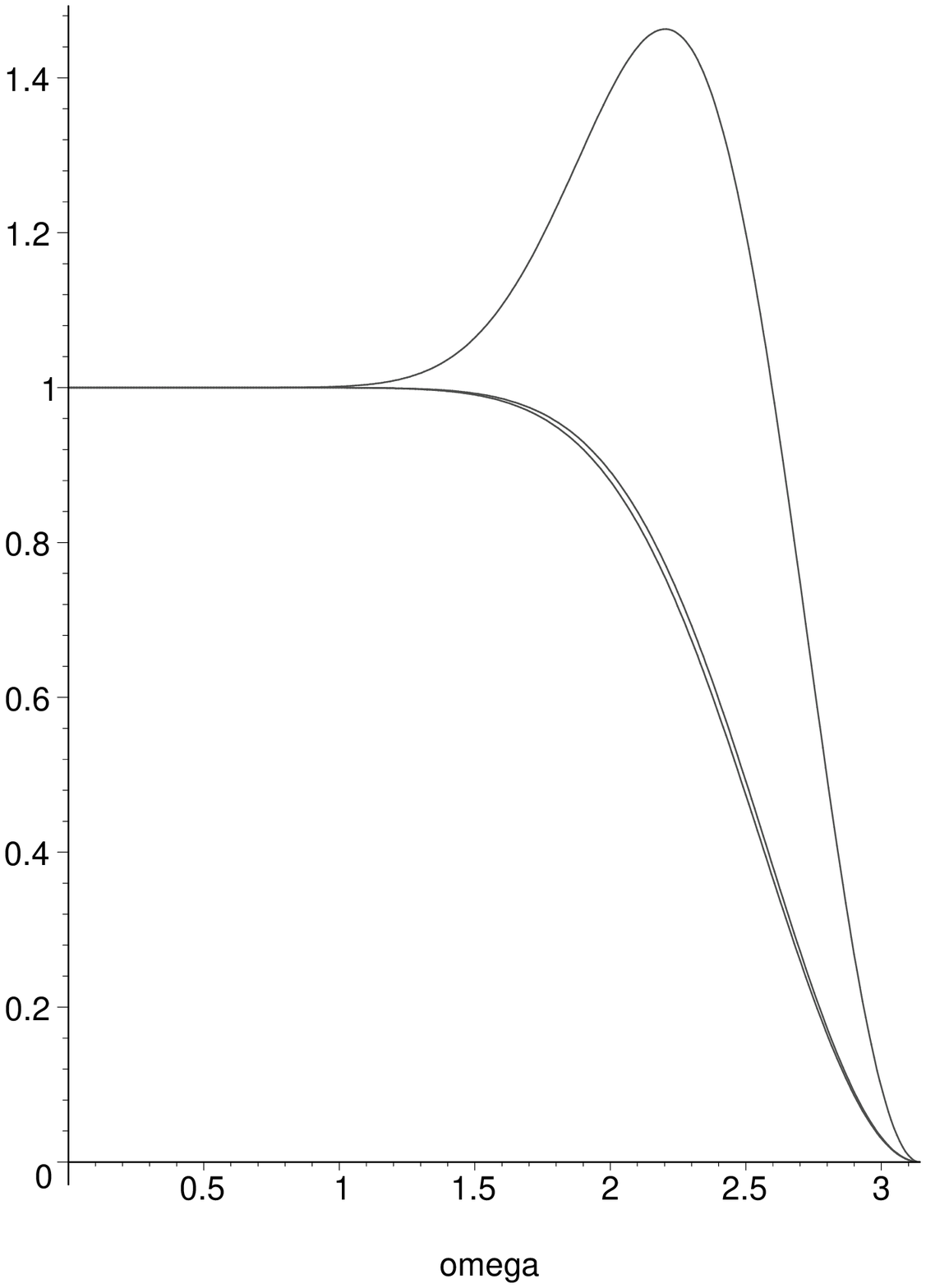}}
\centerline{\smc Figure 1.}

\bigskip
The case we treated above: $L = 1, M = 5, K = 1$ is just within 
Selesnick and Burrus's Region II (see \S 1).
However, ``brute-force'' methods only work in very small cases
in Region II!  For instance, when $L = 0$ it can be seen in several 
different ways that there are $2^M$ complex solutions of the Selesnick-Burrus 
equations.  Thus solving the systems with $L = 0$ 
becomes exponentially more complex as $M$ increases.

\heading \S 3. A Solution Strategy in Region II \endheading

In this section, we will present a strategy for solving the 
Selesnick-Burrus equations in Region II that is much more
efficient than ``brute force'' elimination as in Example 
(2.5).  The idea is to exploit the special structure of the 
Selesnick-Burrus equations as much as possible.
We will also report some results obtained by this strategy.

First, following Selesnick and Burrus,
we show how to reduce the number of variables from 
$N - 1 = K + L + M$ to $M - L - 1$ and obtain an equivalent
system of equations that
we will call the {\it reduced Selesnick-Burrus system}
for a given collection of parameters $K,L,M$.
The computations involved in these steps are minimal.

The first part of this reduction is to use the simple observation
that the last condition $(1+z^{-1})^K|H(z)$ in Selesnick and Burrus'
formulation implies that the moments $m_0, \ldots, m_{N-1}$ already satisfy
certain linear equations, and hence all of the equations
can be expressed in terms of the moments in the column vector
$\vec{m} = (m_0,\ldots, m_{L+M})^{tr}$.  (As noted before, we will
also normalize $m_0 = 1$.)

To see how this works in detail,
write $H(z) = (1+z^{-1})^K P(z)$, let $\vec{h}$ be the 
column vector $(h[0], h[1],\ldots, h[N-1])^{tr}$, and let 
$\vec{p} = (p[0],p[1],\ldots,p[N-1-K])^{tr}$ be the column vector
of coefficients in $P$.   
Then we have an equation 
$$\vec{h} = T \vec{p},\leqno(3.1)$$
where $T$ is an $N\times (N-K)= (K+L+M+1) \times (L+M+1)$ matrix 
whose rows and columns are shifted
copies of the vector of binomial coefficients ${K\choose j}$, $j = 0,\ldots, K$.

By the definition (2.3) of the moments, we have 
$$\vec{m} = Q\vec{h} = QT\vec{p},\leqno(3.2)$$
where $Q$ is an $(L+M+1)\times (K+L+M+1)$ ``Vandermonde-type'' matrix, whose 
$i$th row is the vector of $i$th powers of the integers $0,1,\ldots, K+L+M$.

Combining (3.1) and (3.2), we have the equality
$$\vec{h} = T(QT)^{-1}\vec{m}$$
Hence, we can express $m_k$ for $k > L+M$ as 
$$m_k = (0,1^k,2^k,\ldots,(L+M)^k)T(QT)^{-1}\vec{m}\leqno(3.3)$$

The second part of this reduction is to use some observations
about the Selesnick-Burrus quadratic equations (2.4a) and (2.4b),
and the affine variety they define over the field $\Bbb C$.
It is well-known that there is nothing special about varieties 
defined by quadrics, but the Selesnick-Burrus equations have a 
{\it very particular} form.  First note that the quadratic Selesnick-Burrus
polynomials do not depend on the parameter $K$.  Let $J_{L,M}$
be the ideal they generate in ${\Bbb C}[m_1,\ldots,m_{2M}]$.
In addition, we have the following observations.

\bigskip
\proclaim{(3.4) Lemma}  Let $V_{L,M} = V(J_{L,M})$ be the affine variety
defined by the Selesnick-Burrus quadrics for a given pair of 
parameters $L,M$.
\roster
\item"a." The Selesnick-Burrus quadrics are homogeneous if 
we assign
$$weight(m_i) = i.$$
\item"b." $V_{L,M}$ contains a rational normal curve passing through 
each of its points.
\item"c." $V_{L,M}$ is a smooth variety in ${\Bbb C}^{2M}$ of
dimension $M - L$. 
\endroster
\endproclaim

\proof  All of these claims are easy consequences of the form of 
the quadrics.  \endproof

\bigskip
In fact, we can see much more about the variety $V_{L,M}$ if 
look at another generating set for the ideal that defines it.  
Before giving the general statement, we again take up the 
case $K = 1, L = 1, M = 5$ considered in Example (2.5).

\bigskip
\noindent{\bf (3.5) Example.} Recall the Selesnick-Burrus
quadrics given in Example (2.5).  If we compute a
{\it lex} Gr\"obner basis for $J_{1,5}$ with $m_{10} > m_9 > \cdots > m_1$ we find:
$$\eqalign{G = 
  \{m_2-m_1^2,\ \ m_3-&m_1^3,\ \ m_4-m_1^4,\cr
   m_6-6 m_1m_5+5m_1^6,\ m_8-8 &m_1m_7 + 112 m_1^3m_5 - 105 m_1^8,\cr
 m_{10}-10 m_1 m_9 + 240 m_1^3 m_7 - 126&m_5^2 + 3780 m_1^5 m_5 + 3675 m_1^{10 }\}, \cr}$$
the Gr\"obner basis $G$ shows a very nice {\it parametrization}
for $V_{1,5}$.  If we let 
$$\eqalign{\varphi : {\Bbb C}^4 &\to {\Bbb C}^{10}\cr
(t,a,b,c) &\mapsto (t,t^2,t^3,t^4,a,6at-5t^6,b,8bt-112at^3+105t^8,c,\cr
            &\qquad\qquad 10 ct - 240 bt^3 + 126 a^2 - 3780 at^5 - 3675 t^{10})\cr}$$
The image of $\varphi$ is precisely $V_{1,5}$. 

We next indicate a connection between the Selesnick-Burrus systems and 
some classical topics in algebraic geometry.  These observations 
are needed here only to verify that the hypotheses of \cite{BEM} 
are satisfied for these systems and can be omitted if the reader is not 
familiar with these concepts.  However,
they motivated a large portion of our work on this problem.  

The related ideal 
$$J'=\langle m_2-m_1^2,m_3-m_1^3,m_4-m_1^4,m_6-6 m_1m_5,m_8-8 m_1m_7,m_{10}-10 m_1 m_9\rangle$$
is equal to the ideal generated by the $2\times 2$ minors of
$$
\pmatrix m_1&m_2&m_3&m_4&m_6 &m_8 &m_{10}\cr
         1  &m_1&m_2&m_3&6m_5&8m_7&10m_9\cr\endpmatrix$$
Hence, $S = V(J')$ is an affine 
{\it 4-fold rational normal scroll} (see, e.g. \cite{H}) -- the 
union of ${\Bbb C}^3$'s spanned by related points on
a rational normal curve of degree 4 and 3 lines.  
Moreover, $V_{1,5}$ is the image of the scroll $S$ under a 
certain upper-triangular automorphism $\alpha$ of ${\Bbb C}^{10}$. 
We also see that the projection of $V_{1,5}$ into the coordinate subspace
${\Bbb C}^8$ with coordinates $m_1,\ldots, m_8$ is itself a rational
scroll of dimension 3.  (It is only the quadratic term $a^2$ in the
last coordinate that keeps $V_{1,5}$ from being a rational scroll itself.)

Similar results hold for all the $V_{L,M}$.
These observations imply that $V_{L,M}$ is a {\it unirational}
variety for all $L, M$.  The additional linear equations define the
affine part of a 0-dimensional linear section of $V_{L,M}$.
Because of this, the Selesnick-Burrus systems fall into the general
context discussed in the paper \cite{BEM}, and we can use the main
theorem there to eliminate variables using resultants (without
using Gr\"obner bases).  We will use this approach in the following.

\bigskip
Our next Lemma establishes an important common feature of all of the 
Selesnick-Burrus systems which we will exploit to define the {\it reduced
Selesnick-Burrus system} for a given set of parameters $K,L,M$ with $M > L$.

\bigskip
\proclaim{(3.6) Lemma}  Assume $M > L$. 
The Selesnick-Burrus quadrics (2.4a) and (2.4b)
imply that $m_k = m_1^k$ for all $k$, $1 \le k \le 2L + 2$.
\endproclaim

\proof The proof is by induction on $L$, the base case being $L = 1$.
In that case, we have from (2.4a) with $M = 1,2$, and $m_0 = 1$:
$$2m_1^2 - 2m_2, \quad 6m_2^2 + 2 m_4 - 8 m_1 m_3\leqno(3.7)$$
From (2.4b) with $L = 1$:
$$-m_3 + m_1 m_2\leqno(3.8)$$
The equation $m_2 = m_1^2$ follows directly from the first equation in 
(3.7).  Substituting in (3.8), we have $m_3 = m_1^3$.  Then substituting
in the second equation in (3.8), we have $m_4 = m_1^4$.

The induction step is similar.  Assume we have shown that the quadrics
(2.4a) with $1 \le i \le M$ and (2.4b) $1 \le j \le L$ imply $m_k = m_1^k$
for all $1 \le k \le 2L + 2$.  Consider these quadrics, plus $(2.4a)$ with
$i = M+1$ and (2.4b) with $j = L+1$.  By the induction hypothesis, 
we substitute  $m_k = m_1^k$ for $1 \le k \le 2L + 2$.  Then substituting into
(2.4b) with $j = L+1$, we have 
$$\eqalign{0 &= \left(\sum_{\ell=0}^L \left(1-{2\ell\over 2L+3}\right) 
{2L+3\choose \ell} (-1)^\ell \right) m_1^{2L+3}\cr
     &\qquad  + (-1)^{L+1} \left(1-{2L+2\over 2L+3}\right){2L+3\choose L+1} m_{2L+3}\cr}.$$
This implies $m_{2L + 3} = m_1^{2L+3}$ because, applying some standard
binomial coefficient identities,
$$\sum_{\ell=0}^{L+1} \left(1-{2\ell\over 2L+3}\right) {2L+3\choose \ell} (-1)^\ell 
= \sum_{\ell=0}^{L+1} (-1)^\ell \left({2L+2\choose \ell} - {2L+2\choose\ell-1}\right)= 0.$$
Then, we substitute $m_k = m_1^k$ for $k = 1,\ldots, 2L+3$ into (2.4a) with $i = 2L+4$
to deduce $m_{2L+4} = m_1^{2L+4}$.  \endproof

\bigskip
\noindent
{\bf (3.9) Definition.}  The {\it reduced Selesnick-Burrus system} 
for given parameters $K,L,M$ is the system of equations obtained from 
the full Selesnick-Burrus system of $N-1 = K + L + M$ equations
in $N-1$ variables $m_1, \ldots m_{N-1}$ as follows.
\roster
\item First, substitute in the equations (2.4a) for $i \ge L + 2$,
for all the moments $m_k$ for $k > L + M$ from equations (3.3) above.
\item Write $m_1 = t$ and substitute $m_k = t^k$ for all $1 \le k \le 2L + 2$
in these equations.  Also set $m_0 = 1$.
\endroster
The result is a system of $M - L - 1$ equations in the $M - L - 1$
variables 
$$t,m_{2L+3},\ldots,m_{L+M}.$$
The quadrics (2.4a) with $i \le L+1$ and all of the
quadrics (2.4b) are discarded since they have been used to derive the equations
$m_k = t^k$.
\bigskip

The Gr\"obner basis computation we used in Example (2.5) and substitution of
the parametrization of the variety defined by the Selesnick-Burrus
quadrics does the same sort of elimination of variables as given in 
part 2 of the reduction described here (and more).  Note that the linear 
equations we discussed above, for instance in Example (2.5), have been subsumed 
in the equations (3.3).  We have eliminated the higher moments $m_k$, $k > L+M$ 
using them, so they do not appear explicitly in the reduced system.  
The parameter $K$ enters only in the form of the
$T$ matrix in (3.3).  Changing $K$ changes the coefficients of the equations
but not their Newton polytopes or the number of solutions (provided 
$K \ge 1$).

It will be most useful to view the polynomials in the reduced system as 
polynomials in the moments $m_{2L+3},\ldots,m_{L+M}$, whose
coefficients are polynomials in $t$.  For the Region I cases considered
by Selesnick and Burrus, these polynomials are {\it linear} in $m_{2L+3},\ldots,m_{L+M}$,
and this is what allows the use of purely linear algebra techniques to eliminate
and obtain a univariate polynomial in $t$.  

In fact the Region II cases are characterized by the fact that the 
reduced system still has non-linear terms in $m_{2L+3},\ldots,m_{L+M}$.
The precise form of the reduced system is determined by ``how far down into'' 
Region II we are from the boundary. That is, for $L$ sufficiently large, 
all the cases along the ``diagonals'' defined by $M = 2L + q$, for fixed $q \ge 3$
will have a similar shape.  (There are also ``special cases'' along early portions 
of lower diagonals $M = 2L + q$ with $q \ge 5$.  These are different from the stable 
form because the nonlinear terms are different.)

\bigskip
\noindent
{\bf (3.10) Strategy.}  To study the Selesnick-Burrus equations for 
cases in Region II, we propose the following strategy.
\roster
\item Form the reduced system as in (3.9), and view it as a system of 
$M - L - 1$ linear and quadratic equations in the $M - L - 2$ variables 
$m_{2L+3},\ldots,m_{L+M}$, with the variable $t$ ``hidden in the 
coefficients.''
\item Use the linear equations in the reduced system to solve 
for a subset of the remaining higher moments in terms of the lower moments,
and substitute into the quadratic equations.
\item Use an appropriate formulation for multipolynomial resultants to 
eliminate the remaining undetermined moments and produce a univariate 
polynomial in $t$.  
\endroster

\bigskip
In order to compute examples, we have used several different resultant
formulations.  For instance, in \S 5 below, we will see that the cases with 
$M = 2L + 3$ can be handled by using the multipolynomial resultant of a general system
of $L + 1$ homogeneous linear equations and $1$ homogeneous 
quadratic equation in $L+2$ variables.  
This resultant is denoted by $Res_{1,\ldots,1,2}$ in \cite{CLO}, Chapter 3.

Mixed sparse resultants (see \cite{CE}, \cite{S}), 
Dixon (or B\'ezout) resultants (see \cite{KSY}),
and even the naive approach of iterated pairwise Sylvester resultants 
all work reasonably well on the smaller examples.   
Dixon resultants seem to be far superior for the larger cases.  
In almost all cases, some care is needed to eliminate extraneous factors
in the computed polynomial in $t$.  One useful criterion here is the 
fact mentioned above in (2.6) that the Selesnick-Burrus system is invariant
under time-reversal.  Thus correct univariate polynomial in $t$ must
be invariant under $t \mapsto (K + L + M) - t$.
This strategy is particularly well adapted for the problem of determining
the number of complex solutions of the design equations as a function
of the design parameters $K,L,M$.  In combination with numerical rootfinding
methods, it can also serve as a template for a general solution method
for the Selesnick-Burrus systems.  We illustrate this below.

To indicate the scale of the problems that this strategy allows us
to solve, we provide the following table giving the degree of the 
univariate polynomial in $t$ generating the elimination ideal
of the Selesnick-Burrus system for given $L,M$.  In most cases
the computation was done with $K = 1$ for simplicity, but the 
degree will be the same for all $K \ge 1$.

In this table, the entries along the diagonal $M = 2L + 3$ are the 
first within Region II; the Region I cases with $M < 2L+3$ are not shown.
For purposes of comparison, the entries for $M\le 7$ were also reported by 
Selesnick and Burrus in \cite{SB}.  The entries with $M \ge 8$ and $L > 0$
are new.  Starred entries
were computed by Robert Lewis of Fordham University, using his {\it Fermat} 
system and his routines for Dixon resultants.  The blank entries are somewhat
beyond the scope of current software.  On the other hand, many cases with
$M \ge 15$ would also be tractable by these methods.

$$\vbox{\tabskip=0pt \offinterlineskip
\halign{\strut
\vrule#&\hfil#\ 
&\vrule#&\ \hfil#\ 
&\vrule#&\ \hfil#\ 
&\vrule#&\ \hfil#\ 
&\vrule#&\ \hfil#\ 
&\vrule#&\ \hfil#\ 
&\vrule#&\ \hfil#\ 
&\vrule#&\ \hfil#\ 
&\vrule#\cr
\multispan{14}\hrulefill\cr
&M/L&&0&&1&&2&&3&&4&&5&&\cr\multispan{14}\hrulefill\cr
& 5&&32&& 16&&  &&  &&  &&  &\cr\multispan{14}\hrulefill\cr
& 6&&64&& 26&&  &&  &&  &&  &\cr\multispan{14}\hrulefill\cr
& 7&&128&& 48&&24&&  &&  &&  &\cr\multispan{14}\hrulefill\cr
& 8&&256&& 78&&38&&  &&  &&  &\cr\multispan{14}\hrulefill\cr
& 9&&512&&152&&66&&32&&  &&  &\cr\multispan{14}\hrulefill\cr
&10&&1024&&278&&112&&50&&  &&  &\cr\multispan{14}\hrulefill\cr
&11&&2048&&512*&&192*&&86&&40&&   &\cr\multispan{14}\hrulefill\cr
&12&&4096&&944*&&358*&&142&&62&&  &\cr\multispan{14}\hrulefill\cr
&13&&8192&&   &&572*&&240*&&106&&48&\cr\multispan{14}\hrulefill\cr
&14&&16384&&  &&1020*&&402*&&174*&&74&\cr\multispan{14}\hrulefill\cr
}}
$$
\centerline{\smc Figure 2.}

\bigskip
We will now present an outline of the resultant computation
for the case $K = 2, L = 2, M = 10$ and show how the methods
described in \cite{BEM} and \cite{M} can be used to derive
all the real solutions.  The reduced Selesnick-Burrus system in
this case is a system of $M - L - 1 = 7$ equations in the 
$7$ variables $t = m_1$, and the $m_j$, $j = 7, \ldots, 12$.  
We will begin by using the resultant to eliminate $m_j$, $j = 7,\ldots,12$ and 
yield a univariate polynomial in $t$ satisfied by all the 
solutions.  This is done by ``hiding the variable $t$ in 
the coefficients'' of the system as described, for instance,
in \cite{M}.

For simplicity,
we will write $m_7 = x$, $m_8 = y$, $m_9 = z$, $m_{10} = u$, 
$m_{11} = v$, $m_{12} = w$, and denote the $j$th equation
by $a_j(x,y,z,u,v,w) = 0$.  The first three equations are
$$\eqalign{
0 = a_1(x,y,z,u,v,w) &= 7t^8 + y - 8tx\cr
0 = a_2(x,y,z,u,v,w) &= -84t^{10} - u + 10tx - 45 yt^2 + 120 xt^3\cr
0 = a_3(x,y,z,u,v,w) &= 462 t^{12} + w - 12tv + 66t^2u - 220zt^3 + 495 yt^4 - 792t^5x\cr}$$
The remaining four equations are significantly more complicated and will
be omitted here.  (The complete computation is available as a Maple 8
worksheet from the author's homepage by downloading
\medskip
\centerline{\tt mathcs.holycross.edu/$\sim$little/SB2210.mws}
\medskip
\noindent
To run this and other examples, the procedures in the file
\medskip
\centerline{\tt mathcs.holycross.edu/$\sim$little/CompFileLatest.map}
\medskip
\noindent
should also be downloaded.)

The Dixon resultant computation proceeds as follows.
We introduce a second set of variables $X,Y,Z,U,V,W$ and compute
the $7\times 7$ determinant $\Delta$ whose $j$th row is the transpose of
$$
\pmatrix a_j(x,y,z,u,v,w) \cr
{a_j(X,y,z,u,v,w) - a_j(x,y,z,u,v,w)\over X - x}\cr
{a_j(X,Y,z,u,v,w) - a_j(X,y,z,u,v,w)\over Y - y}\cr
{a_j(X,Y,Z,u,v,w) - a_j(X,Y,z,u,v,w)\over Z - z}\cr
{a_j(X,Y,Z,U,v,w) - a_j(X,Y,Z,u,v,w)\over U - u}\cr
{a_j(X,Y,Z,U,V,w) - a_j(X,Y,Z,U,v,w)\over V - v}\cr
{a_j(X,Y,Z,U,V,W) - a_j(X,Y,Z,U,V,w)\over W - w}\cr\endpmatrix
$$
The expanded form of the determinant can be written as a matrix product
$\Delta = R\cdot M\cdot C$, where $R$ is a 44-component row vector containing 
monomials in $x,y,z,u,v,w$, $M$ is a $44\times 36$ matrix
whose entries are polynomials in $t$, and $C$ is 36-component
column vector whose entries are monomials in $X,Y,Z,U,V,W$.  
The rank of the matrix $M$ in this case is 24.

By the main result of \cite{BEM}, 
any $24\times 24$ submatrix $M'$ of $M$ of rank 24 has
determinant equal to a multiple of the resultant of the system. 
For a particular choice of maximal rank submatrix, we computed 
and factored the determinant yielding a reducible polynomial with 
one factor of degree 112 in $t$ and other factors of smaller degrees.  
The factor of degree 112 is the resultant; the others are extraneous
factors that depend on the choice of the submatrix $M'$.

Using Maple's {\tt fsolve} routine, 12 approximate real roots were 
determined, $t \doteq$
$$\eqalign{
 &.021826159039817\cdots,\cr
1&.14111245031295\cdots,\cr
2&.46849175059426\cdots,\cr
4&.77577862421111\cdots,\cr
5&.42248255383217\cdots,\cr
6&.63285847397435\cdots\cr}$$
and six additional roots obtained from these by time
reversal -- $t \mapsto 14 - t$
(note that $K + L + M = 2 + 2 + 10 = 14$).  In this computation,
a 170 decimal digit floating-point number system was necessary to obtain
accurate results.  The use of the moment variables in the 
Selesnick-Burrus formulation simplifies 
the form of the equations immensely and makes the symbolic
approach we have used feasible.  But it also imposes a severe
numerical conditioning penalty in return.

To determine the other components of the solution, we use the
form of the row monomial vector $R$ in the the equation 
$\Delta = R\cdot M\cdot C$ above.  The entries of $R$ corresponding
to the rows of the maximal rank $24\times 24$ submatrix $M'$ 
contain the six monomials $x,y,z,u,v,w$. 
Substituting each of the $t$ values above in turn, 
the vector in the kernel of $(M')^{tr}$ with 
first component equal to $1$ has 6 components equal to 
the $x,y,z,u,v,w$ values in the corresponding solution of the
system.  We then determine the values of the filter coefficients
from the moments from (3.2) and (3.3) above.  

The square magnitude responses of the 6 real filters found above
are shown in Figure 3.  
Of these, four are apparently monotone decreasing, one has 
a maximum, and one has a minimum and a maximum.  The four
monotone filters come from the $t$-values closest to the
center value $t= {K+L+M\over 2} = 7$.

\bigskip
\epsfxsize 2.0truein
\centerline{\epsfbox{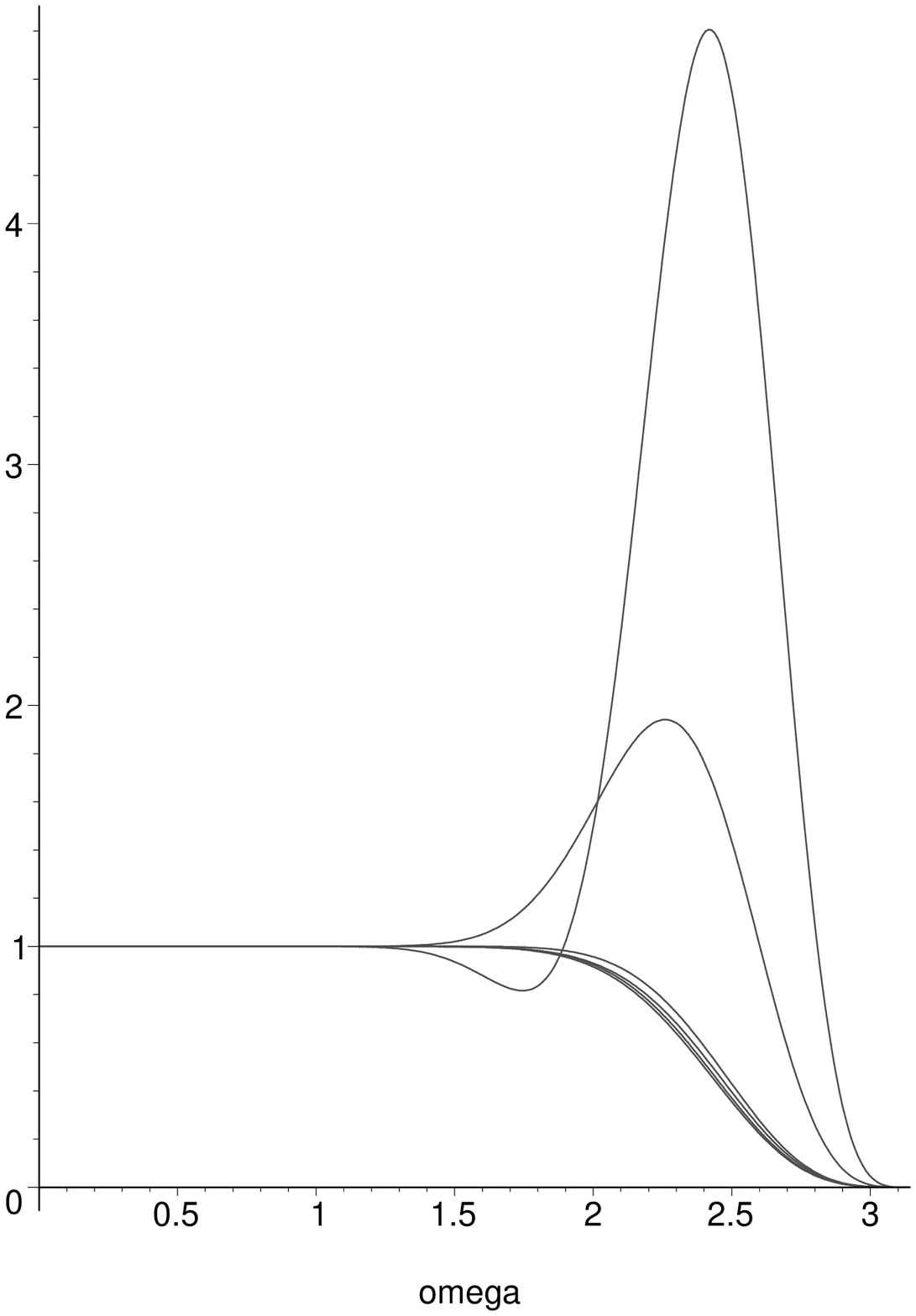}}
\centerline{\smc Figure 3.}
\bigskip

Timings for this computation are as follows (all done in Maple 8
on a SunBlade 100 workstation with a 500 MHz UltraSPARC processor
and 256MB of RAM, running Solaris).  The symbolic
part of the computation (the computation of the Dixon resultant and
factoring the univariate polynomial) takes approximately 320 seconds.
(There is a certain amount of randomness built into the choice
of the maximal rank submatrix $M'$, however, and the time can 
vary depending on which submatrix is used.)
The numerical part (the rootfinding steps) can be done quickly (i.e in less
CPU time than the symbolic computation, even with the high-precision
arithmetic) with an {\it ad hoc} ``by-hand''
search for the real roots in the interval $[0,7]$ and a fast iterative 
method like Newton-Raphson.  (With the ``brute-force'' application of 
Maple's {\tt fsolve} command described above, and illustrated in the 
worksheet mentioned before, the numerical part of the
computation takes much longer, of course -- about 8200 seconds, including
the plotting of the square magnitude response curves.)

We have used similar numerical computations to
solve the reduced system and determine the filter coefficients of the 
real solutions in many of the cases reported above in Figure 2.
As is indicated by this example, we note that the degrees give only one 
measure of the complexity of these computations.  

In a companion paper \cite{LL}, we discuss some properties of the filters
obtained by these computations in more detail.
In the next sections here we will focus instead on some of the
patterns that seem to appear when the table in Figure 2 is examined carefully.

\heading \S 4. A Technical Interlude\endheading

In this section we will prove a number of technical lemmas on the 
Smith normal form of certain matrices that appear when
the linear equations in the reduced Selesnick-Burrus
systems (3.9) are reformulated in a particularly useful way.  
For simplicity, we will describe the general form of these matrices
in this section in the abstract, so to speak; we will delay showing
how the Selesnick-Burrus equations fit these patterns until \S 5
and \S 6.  

We will need the following notation.

\bigskip
\noindent
{\bf Notation.}  Let $j,K,\ell$ denote nonnegative integers, and $t$ an 
indeterminate.  All vectors are infinite, indexed by the nonnegative integers,
${\Bbb Z}_{\ge 0}$.
\bigskip
\item{a.}  We will write $\Delta^j$ for the {\it vector of coefficients in
the $j$th forward difference operator}, each entry divided by $j!$, 
``padded'' with additional zero entries on the right:
$$\Delta^j = {1\over j!} \left((-1)^j {j\choose 0}, (-1)^{j-1}{j\choose 1}, (-1)^{j-2}{j\choose 2},\ldots, 
{j\choose j}, 0, \ldots\ \right).$$
The indices of the nonzero entries shown run from $0$ to $j$.
\item{b.} Similarly, we will write $\Delta^j_\ell$ for right shift by $\ell$ of the 
vector above, so the ${1\over j!} (-1)^j {j\choose 0}$ occurs in position $\ell$, and zeroes
appear in locations $0$ through $\ell - 1$.  
\item{c.}  We will write 
$$D^j_K = {1\over 2^K} \sum_{\ell = 0}^K {K \choose \ell} \Delta^j_\ell.$$
The vector $D^j_K$ can also be viewed as the padded vector of coefficients
of a difference operator.
\item{d.} We will write $(i - t)^\ell$ for the vector with entries
$$((0-t)^\ell, (1-t)^\ell, (2-t)^\ell, \ldots )$$
\item{e.}  We will use the shorthand
$$[j,\ell;K] = \langle D^j_K, (i-t)^\ell\rangle,$$
where $\langle, \rangle$ is the formal dot product on vectors indexed
by ${\Bbb Z}_{\ge 0}$.  Note that all of the vectors $D^j_K$
we consider have only a finite number of nonzero terms, so convergence is automatic.
The sum is the value at $i = 0$ of the result of applying the 
operator $D^j_K$ to the function of the discrete variable $i$ given 
by $(i - t)^\ell$ .  This is a polynomial of degree $\ell - j$ in $t$ if $\ell \ge j$, 
and equals zero otherwise because all $j$th differences of a polynomial of 
degree $< j$ in $i$ vanish.
\item{f.}  An expression of the form $[j,\ell;K](a)$ will denote the value 
obtained by substituting $t = a$ in the polynomial $[j,\ell;K]$.

\bigskip
\proclaim{(4.1) Lemma} The $[j,\ell; K]$ polynomials have the following
properties.
\item{a.} (``reflection identity'') Up to a sign, $[j,\ell;K]$ is symmetric
about $t = {j+K\over 2}$:
$$[j,\ell;K](j+K-t) = (-1)^{j+\ell}[j,\ell;K](t)$$
We call $t = {j + K \over 2}$ the {\bf center value} of $[j,\ell;K]$.
\item{b.} (``center value zero'')  If $\ell$ and $j$ have opposite parity, then 
$$[j,\ell;K]\left({j+K\over 2}\right) = 0.$$
\item{c.} (``boost identity'') $[j,\ell;K]$ satisfies
$$[j,\ell;K](t-1) = [j,\ell;K](t) + (j+1)[j+1,\ell,K](t).$$
\endproclaim

\bigskip
\proof
Part a follows from a direct computation.  Because of the symmetry of the binomial
coefficients in the $\Delta^j_\ell$,  $D^j_K$ is symmetric about ${j+K\over 2}$,
up to the sign $(-1)^j$.  Therefore we have
$$\eqalign{
[j,\ell;K]((j+K)-t) &= \langle D^j_K, (i - ((j+K) - t))^\ell\rangle\cr
                    &= (-1)^\ell \langle D^j_K, (((j+K) - i) - t)^\ell\rangle\cr
                    &= (-1)^{j+\ell}\langle D^j_K, (i - t)^\ell\rangle,\cr
                    &= (-1)^{j+\ell} [j,\ell;K](t).\cr}$$
  
\bigskip
Part b follows immediately from part a.

\bigskip
Part c is shown by another calculation.  In terms of the shift operator $E$, 
we have $D^{j+1}_K = {1\over 2^K (j+1)!} (E+1)^K(E - 1)^{j+1}$, so 
$$\eqalign{(j+1)[j+1,\ell;K](t) &= (j+1)\langle D^{j+1}_K, (i-t)^\ell\rangle \cr
&= {1\over 2^K j!} \langle  (E+1)^K (E - 1)^{j+1},(i-t)^\ell\rangle\cr
&= \langle D^j_K E, (i-t)^\ell\rangle - \langle D^j_K, (i-t)^\ell\rangle\cr
&= [j,\ell;K](t-1) - [j,\ell;K](t).\ \square\cr}$$
\enddemo

\bigskip
The specific matrices that will appear in the analysis of the linear 
equations in the reduced Selesnick-Burrus systems have the following forms
$A(s,m;K)$ and $\tilde{A}(s,m;K)$, for certain
positive integers $s,m$ depending on the flatness parameters $L,M$ from
the filter design problem.   First we 
introduce the matrix $A(s,m;K) = $
$$
     \pmatrix [2s-1,2s;K]&[2s,2s;K]&\ldots &[2s+m-1,2s;K]\cr
              [2s-1,2s+2;K]&[2s,2s+2;K]&\ldots &[2s+m-1,2s+2;K]\cr
              \vdots & \vdots &\ddots & \vdots \cr
              [2s-1,2s+2m;K]&[2s+1,2s+2m;K]&\ldots &[2s+m-1,2s+2m;K]\cr\endpmatrix .\leqno(4.2a)$$
We will write $\delta(s,m;K) = \det A(s,m;K)$.

Similarly,
$\tilde{A}(s,m;K) = $
$$    \pmatrix [2s,2s;K]&[2s+1,2s;K]&\ldots &[2s+m,2s;K]\cr
              [2s,2s+2;K]&[2s+1,2s+2;K]&\ldots &[2s+m,2s+2;K]\cr
              \vdots & \vdots &\ddots & \vdots \cr
              [2s,2s+2m;K]&[2s+1,2s+2m;K]&\ldots &[2s+m,2s+2m;K]\cr\endpmatrix .\leqno(4.2b)$$
We write $\tilde{\delta}(s,m;K) = \det \tilde{A}(s,m;K)$.

\bigskip
For example, with $s = 3$, $m = 1$, and $K = 2$, the matrix  $A(3,1;2)$ is
$$A(3,1;2) = \pmatrix 21-6t&1\cr
-56t^3+588t^2-2212 t+2940&476-224t+28t^2\cr\endpmatrix.$$
The entry in the second row and second column is $[2s,2s+2;K] = [6,8;2]$.

The following observation will simplify our work considerably.

\bigskip
\proclaim{Observation}
Since $[j,\ell;K]$ is zero if $j > \ell$, note that all the entries on
the first row of the matrix $\tilde{A}(s,m;K)$ except the first are zero.
Expanding along the first row, we have 
$$\tilde{\delta}(s,m;K) = \delta(s+1,m-1;K).$$
Therefore, for our purposes it will suffice to study the $\delta(s,m;K)$.
\endproclaim

Our main goal in the remainder of this section is to determine the {\it Smith 
normal form} of the matrices 
$A(s,m;K)$ above, and hence to determine $\delta(s,m;K)$. 
Recall that the Smith normal form of a square matrix $A$ with 
entries in ${\Bbb C}[t]$ is the diagonal matrix obtained by doing 
elementary row and column operations. The diagonal entries satisfy 
the following property for all $n \le rank(A)$:
{\it the product of the first $n$ diagonal entries
is equal to the monic greatest common divisor of all the $n\times n$
minors of $A$}.  The properties of the Smith normal form follow from 
the standard theory of homomorphisms
between modules over a PID such as ${\Bbb C}[t]$ (see for instance \cite{Ja}).

We introduce
the following additional notation to facilitate working column by column
in $A(s,m;K)$.  Note that the entries
in $A(s,m;K)$ all have the form $[j,\ell;K]$ with $2s \le \ell\le 2s + 2m$,
$\ell$ even.
The entries in the first column have $j = 2s - 1$. 
The entries in the second have $j = 2s$, and so forth.
We will write $A^j$ for the column in $A = A(s,m;K)$ in which the entries 
are $[j,\ell;K]$ for $2s \le \ell\le 2s + 2m$, $\ell$ even.
     
 Our first result shows that $\delta(s,m;K)$ is symmetric
about $t = {2s + m - 1 + K\over 2}$, up to a sign.

\bigskip  
\proclaim{(4.3) Lemma}  Let $\delta(s,m;K)$ be as above, and let
$c_A = 2s+m-1+K$ (${c_A\over 2}$ is the center value of the entries of 
the last column in $A(s,m;K)$).  Then
$$\delta(s,m;K)(c-t) = \pm \delta(s,m;K)(t).$$
\endproclaim

\proof
Consider the column $A^{2s+m-1-p}$ for each $0\le p \le m$. The center value for the
entries in this column is $t = {c - p\over 2}$.  By Lemma (4.1), part a, we have 
the corresponding column in $A(s,m;K)(c-t)$ equals
$$\eqalign{
A^{2s+m-1-p}(c-t) &= A^{2s+m-1-p}((c-p)-(t-p))\cr
                  &= (-1)^{2s+m-1-p}A^{2s+m-1-p}(t-p)\cr}$$
Then we apply the ``boost identity'' (Lemma (4.1), part c) repeatedly to deduce
that 
$A^{2s+m-1-p}(t-p)$ equals $A^{2s+m-1-p}(t)$, plus a linear combination of the
terms $A^{2s+m-1-p + q}(t)$, for $1 \le q \le p$.  It follows that the 
the column in $A^{2s+m-1-p}(c-t)$ is in the span of the columns in $A^j(t)$ with 
$2s+m-1-p \le j \le 2s + m - 1$, and hence $\delta(s,m;K)(c-t) = \pm \delta(s,m;K)(t)$.
\endproof

\bigskip
Some factors in $\delta(s,m;K)$ are immediately clear from Lemma (4.1), part b.
If $j$ is odd, then the center value root $t = {j + K\over 2}$ of the entries 
in the column $A^j$ is also a root of $\delta(s,m;K)$.
Moreover, it will will follow from the next lemma
that $\left(t - {j+K\over 2} \right)^e$ with $e > 1$ divides $\delta(s,m;K)$
in some cases.

\bigskip
\proclaim{(4.4) Lemma}  Let $R = [2s-1,2s+m-1] \cap {\Bbb N}$.  If $j \in R$ is 
odd and $j + 2p$ is also in $R$, then 
$$A^j\left({j+K+2p\over 2}\right) \in 
{\text Span}\left\{ A^{j'}\left({j+K+2p\over 2}\right) : j' \in R, j' > j, j'\ {\text even}\right\}.$$
\endproclaim

\proof
Note that ${j+K+2p\over 2}$ is the center value of the entries in the column $A^{j+2p}$.
The proof is a kind of double induction argument -- descending induction on the odd
$j \in R$,  and ascending induction on $p \ge 0$ such that 
$j + 2p \in R$.  In the base case for the outer induction, $j$ is the largest odd 
integer in $R$.  In this case, necessarily, $p = 0$.  But then $t = {j + K \over 2}$ is the 
center value root of the column $A^j$, so the conclusion of the Lemma follows.  
Similarly, if $j$ is any odd integer in $R$ and $p = 0$ we see that 
$A^j\left({j+K\over 2}\right) = 0 \in 
{\text Span}\left\{ A^{j'}\left({j+K+2p\over 2}\right) : j' \in R, j' > j, j'\ {\text even}\right\}$.

For the inductive step, assume that the conclusion of the lemma holds for a given
$j,p$, and also for all odd $\tilde{j} > j$ and all $q$ such that $\tilde{j} + 2q \in R$.
If $j + 2(p+1) \in R$, then we consider $A^j\left({j+K+2(p+1)\over 2}\right)$.  
By the ``boost identity'' from Lemma (4.1), part c, for all even $\ell$, 
$2s \le \ell\le 2s + 2m$, we have $[j,\ell;K]\left({j+K+2(p+1)\over 2}\right)=$ 
$$[j,\ell;K]\left({j+K+2p\over 2}\right) - (j+1)[j+1,\ell;K]\left({j+K+2(p+1)\over 2}\right).$$
Hence 
$$A^j\left({j+K+2(p+1)\over 2}\right) = A^j\left({j+K+2p\over 2}\right) 
- (j+1) A^{j+1}\left({j+K+2(p+1)\over 2}\right).$$
In the second term, $j+1$ is even and $> j$, so we do not need to do anything 
further with that.
We apply the inductive hypothesis to the first term:
$$A^j\left({j+K+2p\over 2}\right) \in 
{\text Span}\{ A^{j'}\left({j+K+2p\over 2}\right) : j' \in R, j' > j, j'\ {\text even}\}.$$
The entries in the $A^{j'}\left({j+K+2p\over 2}\right)$ appearing in the 
linear combination are the 
$[j',\ell;K]\left({j+K+2p\over 2}\right)$.  By the ``boost identity'' from Lemma
(4.1) part c again, we have  $[j',\ell;K]\left({j+K+2p\over 2}\right) =$
$$ [j',\ell;K]\left({j+K+2(p+1)\over 2}\right)
+ (j'+1)[j'+1,\ell;K]\left({j+K+2(p+1)\over 2}\right).$$
Hence $A^{j'}\left({j+K+2p\over 2}\right) \in$
$$ 
{\text Span}\left\{ A^{j'}\left({j+K+2(p+1)\over 2}\right), A^{j'+1}\left({j+K+2(p+1)\over 2}\right)\right\}.$$
In the first vector in this set, $j' > j$ is even and this term matches the conclusion. 
In the second vector, $j'+1 > j$ is odd.  Moreover, since for suitable $q$,
$j+2(p+1) = (j'+1)+2q = \tilde{j} + 2q\in R$, we may apply the induction hypothesis to conclude that the second 
vector is also
in ${\text Span}\{ A^{j'}\left({j+K+2(p+1)\over 2}\right) : j' \in R, j' > j, j'\ {\text even}\}$.
\endproof

\bigskip
The main consequence we will draw from this lemma is the following corollary
giving information about the Smith normal form of $A(s,m;K)$ and $\delta(s,m;K)$.

\bigskip
\proclaim{(4.5) Corollary}  Let $p \ge 0$, let $2s - 1 + 2p \in R$, and let
$t = {2s - 1 + 2p + K\over 2}$, the center value 
for the column $A^{2s - 1 + 2p}$.  Then the rank of $A(s,m;K)$ at this $t$ is
at most $m - p$ (i.e. the rank drops by at least $p + 1$ at this $t$).  Hence
$(2t - (2s - 1 + 2p + K))$ divides the last $p + 1$ entries on the diagonal
of the Smith normal form of $A(s,m;K)$, and $(2t - (2s - 1 + 2p + K))^{p+1}$
divides $\delta(s,m;K)$.
\endproclaim

\proof
By standard properties of the Smith normal form, all the claims here follow from
the statement about the rank of $A(s,m;K)$ at $t = {2s - 1 + 2p + K\over 2}$.
That statement follows directly from Lemma (4.4):  At this $t$, the $p + 1$ columns
$A^{2s-1+2q}$, $0\le q \le p$ are all in the span of the remaining
columns of $A(s,m;K)$. \endproof

\bigskip
For future reference, we note that 
Lemma (4.3) (the symmetry of $\delta(s,m;K)$
about $t = c_A = {2s - 1 + m +K\over 2}$ up to sign)
implies the existence of additional roots of $\delta(s,m;K)$
greater than $c_A$.

The foregoing establishes lower bounds on the multiplicities of 
the roots of $\delta(s,m;K)$ at the 
center value roots of the columns $A^j$ for odd $j$.
We next show that there are also roots of $\delta(s,m;K)$ at the center 
value $t$-values of the columns $A^j$ for even $j$.

\bigskip
\proclaim{(4.6) Lemma}  Let $2s + 2p \in R$ and consider the center value  
$t = {2s + 2p + K\over 2}$ for the column $A^{2s + 2p}$.  If $j$ is odd and
$j < 2s + 2p$, then 
$$A^j\left({2s+2p+K\over 2}\right) \in 
{\text Span}\left\{ A^{j'}\left({2s+2p+K\over 2}\right) : j < j' \le 2s+2p, j'\ {\text even}\right\}.$$
\endproclaim

\proof
The proof is similar to the proof of Lemma (4.4) except that now we will 
proceed by ascending induction on $p$, and descending induction on $j$
such that $j < 2s + 2p$.  The base cases are $p = 0$, $j = 2s - 1$, and 
more generally, $p$ arbitrary and $j = 2s + 2p - 1$.  By the ``boost identity'' ((4.1) 
part c),
$[2s+2p-1,\ell;K]\left({2s + 2p + K\over 2}\right) = $
$$[2s+2p-1,\ell;K]\left({2s + 2p + K\over 2} - 1\right)
-(2s+2p)[2s+2p,\ell;K]\left({2s + 2p + K\over 2}\right)\leqno(4.7)$$
Next, apply the ``reflection identity'' ((4.1) part a) to the first term on the right.
The center value for $A^{2s + 2p - 1}$ is ${2s + 2p - 1 + K \over 2}$, so 
$$2s + 2p - 1 + K - \left({2s + 2p + K\over 2} - 1\right) = {2s + 2p + K\over 2}.$$
Hence, since $\ell$ is even and $2s+2p-1$ is odd, 
$$[2s+2p-1,\ell;K]\left({2s + 2p + K\over 2} - 1\right) = 
-[2s+2p-1,\ell;K]\left({2s + 2p + K\over 2}\right)\leqno(4.8)$$
Combining (4.7) and (4.8) for all even $\ell$, $2s \le \ell \le 2s + 2m$, we have
$$A^{2s+2p-1}\left({2s+2p+K\over 2}\right) \in {\text Span}\left\{ A^{2s+2p}\left({2s+2p+K\over 2}\right)\right\}.$$
So the conclusion of the Lemma holds in these cases.

For the inductive step, assume that the conclusion of the Lemma holds for 
for all odd integers $\tilde{j}$ between $j + 2$ and $2s + 2p$
with the current $p$, and for all $\tilde{p} < p$.
Consider the entries $[j,\ell,K]$ in $A^j$.  By the ``boost identity'' ((4.1) part c),
$[j,\ell,K]\left({2s+2p+K\over 2}\right) = $
$$[j,\ell,K]\left({2s+2p+K\over 2}-1\right)-(j+1)[j+1,\ell,K]\left({2s+2p+K\over 2}\right).$$
Hence
$$A^j\left({2s+2p+K\over 2}\right) \in {\text Span}\left\{ A^j\left({2s+2p+K\over 2}-1\right),
A^{j+1}\left({2s+2p+K\over 2}\right)\right\}.$$
In the second vector on the right, $j+1 > j$ is even so this term matches the conclusion of the 
Lemma.  In the first vector, ${2s+2p+K\over 2}-1 = {2s + 2(p-1) + K\over 2}$ is 
the center value for the column $A^{2s + 2(p-1)}$, and $j < 2s + 2(p-1)$.
By the induction hypothesis, $A^{j}\left({2s + 2(p-1) + K\over 2}\right) \in$
$${\text Span}\left\{ A^{j'}\left({2s+2(p-1)+K\over 2}\right) : j < j' \le 2s+2(p-1), j'\ {\text even}\right\}.$$
But for each entry in one of these $A^{j'}$, we can apply the ``boost identity'' 
again:  $[j',\ell;K]\left({2s+2(p-1)+K\over 2}\right) = $
$$[j',\ell;K]\left({2s+2p+K\over 2}\right) + (j'+1)[j'+1,\ell;K]\left({2s+2p+K\over 2}\right).$$
Hence
$$A^{j'}\left({2s+2(p-1)+K\over 2}\right) \in {\text Span}\left\{ A^{j'}\left({2s+2p+K\over 2}\right),
A^{j'+1}\left({2s+2p+K\over 2}\right)\right\}.$$
The first vector in the spanning set matches the conclusion of the Lemma since $j'$ is even.
In the second term, $j'+1 > j$ is odd.  Hence that vector can be written as a linear
combination as in the conclusion of the Lemma by the induction hypothesis.  \endproof

\bigskip
Here too, the main consequence we will draw from this lemma is a corollary
giving information about the Smith normal form of $A(s,m;K)$ and $\delta(s,m;K)$.

\bigskip
\proclaim{(4.9) Corollary}  Let $p \ge 0$, let $2s + 2p \in R$, and let
$t = {2s + 2p + K\over 2}$, the center value 
for the column $A^{2s + 2p}$.  Then the rank of $A(s,m;K)$ at this $t$ is
at most $m - p$ (i.e. the rank drops by at least $p + 1$ at this $t$).  Hence
$(2t - (2s + 2p + K))$ divides the last $p + 1$ entries on the diagonal
of the Smith normal form of $A(s,m;K)$, and $(2t - (2s + 2p + K))^{p+1}$
divides $\delta(s,m;K)$.
\endproclaim

\proof
As in the proof of Corollary (4.5), all the claims here follow from
the statement about the rank of $A(s,m;K)$ at $t = {2s + 2p + K\over 2}$.
That statement follows directly from Lemma (4.6):  At this $t$, the $p + 1$ 
columns $A^{2s-1+2q}$, $0\le q \le p$ are all in the span of the remaining
columns of $A(s,m;K)$. \endproof

\bigskip
As in the case of the center value zeroes from Corollary (4.5), 
Lemma (4.3) (the symmetry, up to a sign, of $\delta(s,m;K)$
under $t \mapsto c_A - t$, where $c_A = 2s - 1 + m + K$)
implies the existence of a second, symmetrically located collection of roots 
of $\delta(s,m;K)$ greater than $c_A$.  We are now ready for the major 
result of this section.

\bigskip
\proclaim{(4.10) Theorem} Let $c_A = 2s - 1 + m + K$ as above.  
The determinant $\delta(s,m;K)$ can be
written in the form:
$$\delta(s,m;K) = a \prod_{i=2s-1+K}^{2s-1+K+2m} (2t - i)^{\lfloor {m - |c_A - i|\over 2} \rfloor + 1}
\leqno(4.11)$$
for some constant $a$.  In the Smith normal form of $A(s,m;K)$, the $(m+1,m+1)$ entry is 
(a constant times) the product $\prod_{i=2s-1+K}^{2s-1+K+2m} (2t - i)$ (one factor for each root), 
the $(m,m)$ entry is a divisor of this polynomial whose roots are the roots 
of $\delta(s,m;K)$
of multiplicity $\ge 2$, and so forth.
\endproclaim

\bigskip
The $|c_A - i|$ in the exponent ensures the symmetry of the exponents in this expansion about
$c_A$.  To make this somewhat intricate statement more intelligible, before proceeding to the proof, 
we give a small example.  Consider the $4\times 4$ matrix $A(2,3;2)$, which has the 
shape
$$A(2,3;2) = \pmatrix [1] & [0] & 0 & 0 \cr
                      [3] & [2] & [1] & [0] \cr
                      [5] & [4] & [3] & [2] \cr
                      [7] & [6] & [5] & [4] \cr \endpmatrix.$$
Here and in the rest of this article we will use the standing
notational convention:

\bigskip
\noindent
{\bf Notation.}  $[d]$ is shorthand for a polynomial of degree exactly $d$ in $t$. 

\bigskip
\noindent
For instance, the entry $[3]$ in the first column is the polynomial
$[3,6;2] = 425 - 420t + 150t^2 - 20t^3$. The entries marked 0 are
actual zeroes.

According to the 
formula in the statement of the Theorem,  the center value
of the $4$th column is $t = {c_A\over 2}$. where 
$c_A = 2s - 1 + m + K = 2\cdot 2 - 1 + 3 + 2 = 8$.
The set of roots is symmetric about $t = 4$.  The ``predicted'' value for $\delta(2,3;2)$ is 
$$\delta(2,3;2) = a (2t - 5)(2t - 6)(2t - 7)^2(2t - 8)^2(2t - 9)^2(2t - 10)(2t - 11)$$
for some constant $a$.
Using the computer algebra system Maple, we find 
$$\delta(2,3;2) = 672\, \left( 2\,t-11 \right)  \left( t-3 \right)  \left( t-5 \right) 
 \left( 2\,t-5 \right)  \left( 2\,t-7 \right) ^{2} \left( 2\,t-9
 \right) ^{2} \left( t-4 \right) ^{2},$$
and the Smith normal form of $A(2,3;2)$ is:
$$\pmatrix 1&0&0&0\cr
           0&1&0&0\cr
           0&0&p_3(t) &0\cr 0&0&0& p_7(t) \cr \endpmatrix,$$
where 
$$p_3(t) = {t}^{3}-12\,{t}^{2}+{\frac {191}{4}}\,t-63 = {1\over 8}(2t - 7)(2t - 8)(2t - 9)$$
and 
$$p_7(t) = {t}^{7}-28\,{t}^{6}+{\frac {665}{2}}\,{t}^{5}-2170\,{t}^{4}+{\frac {
134449}{16}}\,{t}^{3}-{\frac {77203}{4}}\,{t}^{2}+{\frac {389415}{16}}
\,t-{\frac {51975}{4}}.$$
$p_7(t)$ is the monic polynomial with roots $t = 5/2, 3, 7/2, 4, 9/2, 5, 11/2$ (all
multiplicity 1).  We now proceed to the proof of Theorem (4.10).

\proof  It follows from Corollaries (4.5) and (4.9) that the 
product in equation (4.11) divides $\delta(s,m;K)$.  If we knew that $\delta(s,m;K)$
had the form given in (4.11), then the claims about the Smith normal form of 
$A(s,m;K)$ would also follow from these Corollaries.  Hence, to prove the Theorem
it suffices to prove that the degree of $\delta(s,m;K)$ equals the degree of the 
product in (4.11) in $t$.  To compute the degree of $\delta(s,m;K)$, recall the form 
of the matrix $A(s,m;K)$ given in (4.2a).  We have
$$A(s,m;K) = \pmatrix [1] & [0] &  0  &  0  & \cdots & 0\cr
                      [3] & [2] & [1] & [0] & \cdots & 0\cr
                      [5] & [4] & [3] & [2] & \cdots & 0\cr
                      \vdots&\vdots &\vdots &\vdots & &\vdots\cr
                    [2m+1]& [2m]& [2m-1]&[2m-2]& \cdots &[m+1]\cr\endpmatrix$$
where, as earlier, $[d]$ denotes a polynomial of degree $d$ in $t$.  The 0 entries are actual
zeroes.  By examining the form of this matrix, it is not difficult to see that because
of the zeroes above the main diagonal, 
every nonzero product of entries, one from each row and one from each column, has
the same total degree as the product of the entries on the main diagonal:
$$1 + 2 + \cdots + (m + 1) = {(m+1)(m+2)\over 2}.$$
Hence the degree of $\delta(s,m;K)$ is no larger than ${(m+1)(m+2)\over 2}$.

But on the other hand, we will see that the product in (4.11) also has degree 
${(m+1)(m+2)\over 2}$.  Hence it follows that $\delta(s,m;K)$ equals the product in (4.11).
To compute the degree of (4.11), we consider the cases $m$ even and $m$ odd separately.
If $m = 2q$ is even, then the central value of the last column of the matrix gives one 
of the central value roots.
The sum of the multiplicities in (4.11) gives
$$2(1 + 1 + 2 + 2 + \cdots + q + q) + q + 1 = (q + 1)(2q + 1) = {(m+1)(m+2)\over 2}$$
Similarly with $m = 2q + 1$ an odd number, the total degree is 
$$2(1 + 1 + 2 + 2 + \cdots + q + q + (q + 1)) + q + 1 = (q + 1)(2q + 3) = {(m+1)(m+2)\over 2},$$
which concludes the proof.
\endproof

\bigskip
By the Observation above concerning the $\tilde{A}(s,m;K)$ matrices, we 
have a parallel formula for $\tilde{\delta}(s,m;K)$.

\bigskip
\proclaim{(4.12) Corollary} Let $\tilde{c} = 2s + m + K$.  
The determinant $\tilde{\delta}(s,m;K)$ can be
written in the form:
$$\tilde{\delta}(s,m;K) = a \prod_{i=2s+1+K}^{2s-1+K+2m} (2t - i)^{\lfloor {m - 1 - |\tilde{c} - i|\over 2} \rfloor + 1}
\leqno(4.13)$$
for some constant $a$.  In the Smith normal form of $\tilde{A}(s,m;K)$, the $(m,m)$ entry is
a constant times 
the product $\prod_{i=2s+1+K}^{2s-1+K+2m} (2t - i)$ (one factor for each root), 
the $(m-1,m-1)$ entry is the divisor of this whose roots are the roots of $\tilde{\delta}(s,m;K)$
of multiplicity $\ge 2$, and so forth.
\endproclaim

\proof
This follows directly from Theorem (4.10), using the relation
$\tilde{\delta}(s,m;K) = \delta(s+1,m-1;K)$.\endproof

\bigskip
Here is an example, showing $\tilde{\delta}(2,3;2)$ for comparison
with $\delta(2,3;2)$ computed earlier.  Using Maple, we have
$$\tilde{\delta}(2,3;2) = 
36\, \left( t-5 \right)  \left( 2\,t-7 \right)  \left( 2\,t-11 \right)  \left( t-4 \right)  \left( -9+2\,t \right) ^{2},$$
which agrees with (4.13) for this $s,m,K$.

\heading \S 5. The $M = 2L + 3$ Diagonal \endheading

In this section we will discuss the Selesnick-Burrus systems for parameters
$L,M$ satisfying $M = 2L + 3$.  In particular, we will prove the following
theorem which explains one pattern that can be seen in the table given in Figure 2.

\bigskip
\proclaim{(5.1) Theorem} In the cases $M = 2L + 3$, $L\ge 0$, 
(the ``corners'' in Region II boundary), for all $K \ge 1$, the univariate
polynomial in $t$ in the elimination ideal of the Selesnick-Burrus
equations obtained via Strategy (3.10) has degree $8L+8$.
\endproclaim

Before giving the details of the proof, we outline the method we will use.
Along the first diagonal in Region II, of the $M - L - 1 = L + 2$ equations in
the reduced Selesnick-Burrus system, $L + 1$ are inhomogeneous
linear equations in $m_{2L+3},\ldots, m_{3L+3}$
whose coefficients are polynomials in $t$.  We will begin
by showing how the coefficient matrix of this linear part of the system
can be rewritten as the matrix $A(L+2,L;K)$ as defined in \S 4, times
a suitable invertible lower-triangular matrix.  The last equation 
(from the flatness condition $F^{(2M)}(0) = F^{(4L+6)}(0) = 0$)
contains the non-linear term $m_{2L+3}^2$, plus linear
terms in $m_{2L+3},\ldots, m_{3L+3}$.  To eliminate to a univariate polynomial 
in $t$, we will use a formula for the multivariable 
resultant $Res_{1,\ldots,1,d}$ from Proposition 5.4.4 of \cite{Jo} (see also 
Exercise 10 of Chapter 3, \S 3 in \cite{CLO}).  (This formula may be proved 
by the basic approach of solving the linear equations for 
$m_{2L+3},\ldots, m_{3L+3}$ in terms of $t$ by Cramer's Rule, then substituting 
into the last equation to obtain a univariate polynomial in $t$.)
We will need to keep careful track of the factorizations of
$\delta(L+2,L;K) = \det A(L+2,L;K)$ from Theorem (4.10).
The steps in this outline will be accomplished in a series of Lemmas.

\bigskip
\proclaim{(5.2) Lemma} The $L+1$ linear equations in the reduced 
Selesnick-Burrus system with $M = 2L+3$ can be rewritten in the form
$$A(L+2,L;K)\cdot {\Cal L}\cdot \vec{m}_r = \vec{b},$$
where $A(L+2,L;K)$ is the matrix defined in (4.2a), ${\Cal L}$ is a
constant lower-triangular matrix with diagonal entries equal to 1,
$\vec{m}_r = (m_{2L+3},\ldots,m_{3L+3})^{tr}$, and 
$\vec{b} = ([2L+4],[2L+6],\ldots,[4L+4])^{tr}$. 
\endproclaim

\proof
Recall the form of the Selesnick-Burrus quadrics from (2.4a):
$$0 = {2j\choose j} m_j^2 + 2\sum_{\ell=0}^{j-1} {2j\choose \ell}(-1)^{j+\ell} m_\ell m_{2j-\ell}\leqno(5.3)$$
The linear equations in the reduced Selesnick-Burrus system
come from these for $j = L+2,\ldots, 2L+2$, via the reduction process described in (3.9). 
We begin by rearranging these equations to the 
following form by separating the terms involving the variables 
$m_{2L+3},\ldots,m_{3L+3}$ from those depending on
the higher moments $m_{3L+4},\ldots,m_{4L+4}$.  We
have
$$W_1 V_1 T(QT)^{-1}\vec{m} + W_2 V_2 T(QT)^{-1}\vec{m} = \vec{b'},\leqno(5.4)$$
where 
\roster
\item the matrix $W_1$ comes from the coefficients of the $m_k$, $2L+3\le k \le 3L+3$ in (5.3):
$$W_1=\pmatrix -{2L+4\choose 1}t & {2L + 4\choose 0} & 0 & 0 & \cdots \cr
         -{2L+6\choose 3}t^3 & {2L + 6\choose 2}t^2 & -{2L+6\choose 1}t & {2L+6\choose 0} & \cdots \cr
          \vdots    &\vdots &\vdots &\vdots & \cr
         -{4L+4\choose 2L+1}t^{2L+1} & {4L+4\choose 2L}t^{2L} & -{4L+4\choose 2L-1}t^{2L-1} & {4L+4\choose 2L-2}t^{2L-2} & \cdots  
          \cr\endpmatrix$$
\item the matrix $W_2$ comes from the coefficients of the $m_k$, $3L+4 \le k \le 4L+4$ in (5.3):
$$W_2 = (-1)^{L} \pmatrix 0 & 0 & \cdots &0\cr
                          \vdots &\vdots & & \vdots\cr
                          {4L+4\choose L} t^L & -{4L+4\choose L-1}t^{L-1} & \cdots & {4L+4\choose 0}\cr \endpmatrix$$ 
\item the matrices $V_1,V_2$ are Vandermonde-type matrices:
$$V_1 = \pmatrix 0 &  1^{2L+3} & \cdots & (3L+3+K)^{2L+3}\cr
               \vdots &\vdots & \ddots & \vdots\cr
               0 & 1^{3L+3} & \cdots & (3L+3+K)^{3L+3}\cr\endpmatrix,$$
and 
$$V_2 = \pmatrix 0 &  1^{3L+4} & \cdots & (3L+3+K)^{3L+4}\cr
               \vdots &\vdots & \ddots & \vdots\cr
               0 & 1^{4L+4} & \cdots & (3L+3+K)^{4L+4}\cr\endpmatrix$$
\item $\vec{m} = (1,t,\ldots,t^{2L+2},m_{2L+3},\ldots,m_{3L+3})^{tr}$
\item $\vec{b'}$ has the same form as $\vec{b}$ in the statement of the Lemma but is not the entire
vector of $t$ terms.  (There are also terms depending only on $t$ that 
come from the matrix product $(W_1V_1+W_2V_2)T(QT)^{-1}\vec{m}$.) 
\item the matrices $Q$ and $T$ are as in the discussion leading up to (3.3).
\endroster 

Since the first $2L+3$ entries of $\vec{m}$ depend only on $t$,
the coefficients of $m_{2L+3}$ through $m_{3L+3}$ in our equations
come from the product 
$$(W_1V_1 + W_2V_2)\cdot{\Cal T}\cdot \vec{m}_r,$$
where ${\Cal T}$ is the submatrix of $T(QT)^{-1}$ containing all the 
entries from the last $L+1$ columns.  The other terms 
in the product $(W_1V_1 + W_2V_2)\cdot T(QT)^{-1}\cdot \vec{m}$
containing only powers of $t$ go into the vector $\vec{b}$, and 
(5.4) becomes
$$(W_1V_1 + W_2V_2)\cdot{\Cal T}\cdot \vec{m}_r = \vec{b}\leqno(5.5)$$

The fact that establishes the connection between these equations and the matrices $A(s,m;K)$
considered in \S 4 is the following observation.  In the matrix ${\Cal T}$, the final column 
is the vector $D^{3L+3}_K$ as in the Notation at the start of \S 4, written as a column.  This 
follows if we think of the columns of $(QT)^{-1}$ as operators acting on the rows of 
$QT$, thought of as power functions of a discrete variable.  Similarly,
the next-to-last column of $T(QT)^{-1}$ is a linear combination $D^{3L+2}_K + \alpha D^{3L+3}_K$
for some constant $\alpha$, and so on.  In general, we have
$${\Cal T} = \left(D^{2L+3}_K | D^{2L+4}_K | \cdots | D^{3L+3}_K\right)\cdot {\Cal L}\leqno(5.6)$$
for a lower-triangular square $(L+1)\times (L+1)$ matrix ${\Cal L}$ with diagonal entries 1. 

To finish the proof of the Lemma, we substitute (5.6) into (5.5) and rearrange the terms again:
$$(W_1V_1 + W_2V_2)\left(D^{2L+3}_K | D^{2L+4}_K | \cdots | D^{3L+3}_K\right)\cdot {\Cal L}\cdot \vec{m}_r + \vec{b}.$$
We have
$${\Cal U}_1 := V_1 \left(D^{2L+3}_K | D^{2L+4}_K | \cdots | D^{3L+3}_K\right) = (D^j_K i^\ell),$$
for $2L+3\le j \le 3L+3$ and $2L+3 \le \ell \le 3L+3$., and
$${\Cal U}_2:= V_2 \left(D^{2L+3}_K | D^{2L+4}_K | \cdots | D^{3L+3}_K\right) = (D^j_K i^\ell),$$
for $2L+3\le j \le 3L+3$ and $3L+4 \le \ell \le 4L+4$.   

Consider the $(I,J)$ entry of the product 
$$(W_1V_1 + W_2V_2)\left(D^{2L+3}_K | D^{2L+4}_K | \cdots | D^{3L+3}_K\right),$$ 
which is the dot product of the $I$th row of $W_1$ with the $J$th column
of ${\Cal U}_1$, plus the dot product of the $I$th row of $W_2$ with the 
$J$th column of ${\Cal U}_2$.  The form of the entries in $W_1$ and $W_2$ on 
the $I$th row is $(-1)^q {2L + 2(I+1)\choose q} t^q$ for $q$ from 
$2I - 1$ down to $0$.  Hence this sum of dot products equals
$$\eqalign{\sum_{q = 0}^{2I-1} (-1)^q {2L + 2(I+1)\choose q} t^q D^J_K i^{2L + 2(I+1) - q} 
&= D^J_K(i - t)^{2L + 2(I+1)}\cr
& = [J,2L+2(I+1);K],\cr}$$
using the notation introduced in \S 4.  As $I$ runs from $1$ to $L + 1$ and
$J$ runs from $2L + 3$ to $3L + 3$, we see that these entries form 
the matrix $A(L+2,L;K)$ as claimed.
\endproof

\bigskip
For a general system 
of $L+1$ linear homogeneous equations and one homogeneous 
quadratic equation in $L + 2$ variables, if the linear equations
are written as ${\Cal A} \vec{x} = 0$, and the quadratic equation
is $Q(\vec{x}) = 0$, then by the result from Proposition 5.4.4 of 
\cite{Jo} mentioned before, the multivariable resultant $Res_{1,\ldots,1,2}$
equals 
$$Q(\delta_1,-\delta_2,\delta_3,\ldots,(-1)^{L+1}\delta_{L+2})\leqno(5.7)$$
where $\delta_I = \det {\Cal A}_I$, and ${\Cal A}_I$ is the $(L+1)\times (L+1)$
submatrix of ${\Cal A}$ obtained by deleting column $I$.  

We apply this to our reduced Selesnick-Burrus system.  
Write the augmented matrix of the linear equations as
$${\Cal A} = (A(L+2,L,K)\cdot {\Cal L} | -\vec{b})$$
where ${\Cal L}$ is the lower triangular matrix 
and $\vec{b}$ is the column vector $([2L+4],[2L+6],\ldots,[4L+4])^{tr}$
from Lemma (5.2).  
Our next Lemma shows that the determinants
of the minors of ${\Cal A}$ have a common factor of degree ${L\choose 2}$.
To prepare for this statement, we introduce the following notation.
Let $\delta$ be the product of the first $L$ diagnonal entries (elementary 
divisors) in the Smith normal form of $A(L+2,L;K)$:
$$\delta = \prod_{i=2L+5+K}^{4L+1+K} (2t - i)^{\lfloor {L - |3L+3+K - i|\over 2} \rfloor}\leqno(5.8)$$
(There is one factor in this product for each of the roots of multiplicity
$\ge 2$ of $\delta(L+2,L;K)$, and the exponents are each 1 less than the
corresponding exponents in $\delta(L+2,L;K)$.)

\bigskip
\proclaim{(5.9) Lemma}  Let ${\Cal A}$ be as in (5.6) and $\delta_i$
be the $i$th minor of ${\Cal A}$ as above.  If $1 \le i \le L+1$, then
$$\delta_i = [4L+3+i]\cdot \delta$$
where $[4L+3+i]$ is some polynomial in $t$ of degree $4L+3+i$, and $\delta$ is 
the product from (5.8). If $i = L+2$, then 
$$\delta_{L+2} = \prod_{i=2L+3+K}^{4L+3+K} (2t - i)^{\lfloor {L - |3L+3+K - i|\over 2} \rfloor + 1} .
= [2L+1]\cdot \delta$$
\endproclaim

\proof
We begin by computing the minor ${\Cal A}_{L+2}$.
Since ${\Cal L}$ is a constant lower-triangular matrix with 
diagonal entries equal to 1, 
${\Cal A}_{L+2} = \delta(L+2,L;K) = \det A(L+2,L;K)$.
We use Theorem (4.10) to compute this.  We have
$c_A = 3L + 3 + K$ and
$$\delta_{L+2} = \delta(L+2,L;K) 
= a \prod_{i=2L+3+K}^{4L+3+K} (2t - i)^{\lfloor {L - |3L+3+K - i|\over 2} \rfloor + 1}
$$
for some constant $a$.  By the properties of the Smith normal form,
we know that at a root $t = t_0$ of multiplicity $r$, the rank of  
$A(L+2,L;K)$ is $L+1-r$, so every $(L+2-r) \times (L+2-r)$
submatrix of $A(L+2,L;K)$ will have zero determinant at $t = t_0$.

Now consider the other minors ${\Cal A}_i$, for $1 \le i \le L+1$, 
and expand the determinant along the column containing the entries
from the vector $\vec{b}$.  Each term in this expansion is the product of an
entry from $\vec{b}$ times the determinant of an $L\times L$
submatrix of $A(L+2,L;K)\cdot {\Cal L}$.  Hence by the statement at
the end of the last paragraph, $\delta_i$ is divisible by $\delta$.
The remaining factor in ${\Cal A}_i$ comes by examining the degrees of the 
entries of the matrix as in the proof of Theorem (4.10).  Note that
if $L = 1$, the starting value of the index $i$ is greater than the final value.
In that case $\delta = 1$.  In all other cases, the degree of $\delta$ is 
${L\choose 2}$ (see the proof of Theorem (4.10)).
\endproof

We now consider the quadratic polynomial in the reduced Selesnick-Burrus
system.  Let $z$ be a homogenizing variable.   Then the homogenized version 
of $Q$, the equation from $F^{(4L+6)}(0) = 0$, has the form
$$[0]m_{2L+3}^2 + [2L + 3] m_{2L+3}z  + \cdots + [L+3] m_{3L+3}z + [4L+6]z^2\leqno(5.10)$$
We analyze the result of substituting the $(-1)^{i+1}\delta_i$ 
into this polynomial as in (5.7). 

\bigskip
\proclaim{(5.11) Lemma}  The resultant of our system has the form
$$[8L+8]\delta^2,$$
where $[8L+8]$ denotes a polynomial of degree $8L + 8$ in $t$, and 
$\delta$ is the product from (5.8).
\endproclaim

\proof
To obtain the resultant of our equations to eliminate $m_{2L+3},\ldots,m_{3L+3}$,
we substitute
$$\eqalign{m_{2L+3} &= \delta_1\cr
                    &\ \vdots \cr
           m_{3L+3} &= \delta_{L+1}\cr
                z   &= \delta_{L+2}\cr}$$
into (5.10) (following equation (5.7) above), and use Lemma (5.9).  We obtain the following 
expression for the resultant:
$$\eqalign{
&[0]([4L+4]\delta)^2 + [2L+3]([4L+4]\delta)([2L+1]\delta) + \cdots + \cr
&\quad [L+3]([5L+4]\delta)([2L+1]\delta) +[4L+6]([2L+2]\delta)^2\cr
& = [8L+8]\delta^2.\cr}\leqno(5.12)$$
\endproof

The factor of degree $8L + 8$ is the univariate polynomial in $t$ that we want,
and this concludes the proof of Theorem (5.1).
The other factor in (5.12) is extraneous in the sense that the $t$ with $\delta(t) = 0$
do not give solutions of the whole Selesnick-Burrus system.  In fact, it
can be seen that the linear equations in the reduced system 
are {\it inconsistent} for those $t$.  
In algebraic geometric terms, the resultant of the homogenized system contains information
about all the solutions of the equations in projective space, including solutions 
``at infinity''.  The common factor $\delta^2$ gives solutions at infinity, and the 
degree in $t$ of the full polynomial in (5.11) is the degree of the projective closure 
of the affine variety
defined by the Selesnick-Burrus quadrics -- the deformed rational scroll as in the discussion 
given in Example (3.5) in the case $L = 1, M = 5$.  In that case there are no 
solutions at infinity (since $\delta = 1$).  However for $L \ge 2$, there are always such 
solutions.  For example with $L = 2$, there are 24 solutions of the Selesnick-Burrus system 
for all $K \ge 1$, but the degree of the variety defined by the quadrics is 26.  The 
factor $\delta^2 = \left[{2\choose 2}\right]^2 = [2]$ accounts for the difference.  
Similarly, with $L = 3$, there are 32 solutions of the Selesnick-Burrus
equations for all $K \ge 1$, but the degree of the variety defined by the quadrics
is 38.  Again, the factor $\delta^2 = \left[{3\choose 2}\right]^2 = [6]$ 
accounts for the difference.

In the companion article \cite{LL}, we will give more details on the structure of
the filters corresponding to the $8L + 8$ solutions of the Selesnick-Burrus
equations for small $L$.  For example, rather extensive calculations
suggest the following conjectures.

\bigskip
\proclaim{(5.12) Conjectures}  Consider the Selesnick-Burrus
equations with $M = 2L + 3$, and $K \ge 1$.
\roster
\item The polynomial of degree $8L+8$ is irreducible over $\Bbb Q$, hence
has 8L+8 distinct solutions in $\Bbb C$.  
\item Of the $8L+8$ solutions, $2(L + 2)$ are real (yielding
$L + 2$ different filters because of the invariance under time reversal).
\item Exactly four of these (2 different filters), those
with $t = m_1$ closest to the center value ${K+L+M\over 2}$, 
yield monotone decreasing square magnitude response.  
\item The other solutions correspond to filters with progressively greater oscillation
and greater maximum ``passband ripple'' as the distance from $t = m_1$ to 
${K+L+M\over 2}$ increases.
\endroster
\endproclaim

The beginnings of this pattern can be seen in Example (3.5), which gives the 
case $L = 1, M = 5$.

\heading \S 6. The $M = 2L + 4$ Diagonal \endheading

In this section we will discuss the Selesnick-Burrus systems with
$M = 2L + 4$, the second diagonal in Region II in the table given in Figure 2.
Our goal is to prove a result parallel to Theorem (5.1) giving
the degree of the univariate polynomial in $t$ whose roots
give the different solutions.

\bigskip
\proclaim{(6.1) Theorem} In the cases $M = 2L + 4$, $L\ge 1$, 
for all $K \ge 1$, the univariate polynomial in $t$ in the elimination 
ideal of the Selesnick-Burrus equations obtained via Strategy (3.10) 
has degree $12L + 14$.
\endproclaim

Our proof will follow the same pattern as the proof of Theorem (5.1).
First, we analyze the form of the equations in these cases.
We rewrite the linear equations in a suitable form making use of 
the results of \S 4.  Then the univariate polynomial is obtained
via an elimination of variables tailored
to the form of these equations.   

We begin by noting that the reduced Selesnick-Burrus system
in these cases has the following form.  The first $L + 1$
equations (from the flatness conditions 
$F^{(2L+4)}(0) = \cdots = F^{(4L+4)}(0)$ 
are linear in the $L + 2$ variables $m_{2L+3}, \ldots, m_{3L+4}$.  
The remaining two equations have nonlinear terms.  The 
condition $F^{(4L+6)}(0) = 0$ gives a reduced equation containing
$m_{2L+3}^2$, plus linear terms in all the variables.  (This is
the same as the last equation in the $M = 2L+3$ cases.)  In addition,
the condition $F^{(4L+8)}(0) = 0$ gives a reduced equation containing
$m_{2L+4}^2$, $m_{2L+3}m_{2L+5}$, terms, plus linear terms.  
Following the
strategy (3.10), we solve the linear equations for 
$L+1$ of the variables in terms of the others, substitute into 
the quadrics, then compute the Sylvester resultant of the 2 quadrics.
(Our approach here is closely related to one way to derive the
multipolynomial resultant for a system of $L+1$ homogeneous linear 
and $2$ homogeneous quadratic equations in $L + 3$ variables:  
$Res_{1,\ldots,1,2,2}$, but it seems to be easier in this case to use an 
{\it ad hoc} approach.) 

We begin with 
the following Lemma describing the linear equations,
Since the precise statement involves some new quantities, we will sketch
the derivation first, then give the formulation of the Lemma we will use.
First, an argument exactly like the proof of Lemma (5.2) shows that 
the linear equations can be rewritten in the form 
$$A \cdot {\Cal L}\cdot \vec{m}_r = \vec{b}$$
where $A$ is the $(L+1)\times (L+2)$ matrix:
$$\pmatrix 
[2L+3,2L+4;K] & [2L+4,2L+4;K] & \cdots & [3L+4,2L+4;K]\cr
[2L+3,2L+6;K] & [2L+4,2L+6;K] & \cdots & [3L+4,2L+6;K]\cr
\vdots        & \vdots        &        & \vdots\cr
[2L+3,4L+4;K] & [2L+4,4L+4;K] & \cdots & [3L+4,4L+4;K]\cr\endpmatrix,$$
${\Cal L}$ is a certain lower-triangular constant matrix with
1's on the main diagonal, 
$\vec{m}_r = (m_{2L+4},\ldots,m_{3L+4})^{tr}$,
and $\vec{b} = ([2L+4],[2L+6],\ldots,[4L+4])^{tr}$.
We will write $\{2L+3,2L+2i+2;K\}$ for the entry in column 1 and 
row $i$ of the matrix $A\cdot {\Cal L}$ (a certain linear combination
of the entries on row $i$ of the matrix $A$).
After we subtract all terms involving $m_{2L+3}$ to the right-hand sides
of the equations, 
we obtain the following result, because the submatrix 
of $A$ consisting of all entries in the last $L+1$ columns
is precisely the matrix $\tilde{A}(L+2,L;K)$ from (4.2b).

\bigskip
\proclaim{(6.2) Lemma} Using the notation introduced above,
the $L+1$ linear equations in the reduced 
Selesnick-Burrus system with $M = 2L+4$ can be rewritten in the form
$$\tilde{A}(L+2,L;K)\cdot {\Cal L}\cdot \vec{m}_r = \vec{b'},$$
where $\vec{b'} = $
$$([2L+4] - \{2L+3,2L+4;K\} m_{2L+3},\ldots,[4L+4] - \{2L+3,4L+4;K\}m_{2L+3})^{tr}.$$
\endproclaim

We can solve the system 
$\tilde{A}(L+2,L;K)\cdot {\Cal L}\cdot \vec{m}_r = \vec{b'}$
for the moments in $m_r$ using Cramer's Rule.  For $1 \le i \le L+1$,
this gives
$$m_{2L+3+i} = {\det A_i\over \tilde{\delta}(L+2,L;K)}\leqno(6.3)$$
where $A_i$ is the matrix obtained from $\tilde{A}(L+2,L;K)\cdot {\Cal L}$
by replacing column $i$ with the vector $\vec{b'}$.  Next, we consider
what happens when we substitute from (6.3) into the first
nonlinear equation (from $F^{(4L+6)}(0) = 0$).  We will show that 
the result is an equation of the form
$$[0]m_{2L+3}^2 + [2L+3]m_{2L+3} + [4L+6] = 0\leqno(6.4)$$
(in other words, the denominators from (6.3) cancel with terms in the numerators
in this equation).
The situation that produces this cancellation is described in 
the following general lemma.

\bigskip
\proclaim{(6.5) Lemma}  Consider a system of equations of the form
$$\eqalign{a_{11}(t)x_1 + a_{12}(t)x_2 + \cdots + a_{1n}x_n &= r_1(t)\cr
           a_{21}(t)x_1 + a_{22}(t)x_2 + \cdots + a_{2n}x_n &= r_2(t)\cr
               \vdots\qquad\qquad & \qquad \vdots\cr
   a_{n-1,1}(t)x_1 + a_{n-1,2}(t)x_2 + \cdots + a_{n-1,n}x_n &= r_{n-1}(t)\cr
   a_{n1}(t)x_1 + a_{n2}(t)x_2 + \cdots + a_{nn}x_n &= r_n(t) + cx_1^2,\cr}$$
where $a_{ij}(t)$ and $r_i(t)$ are in ${\Bbb C}[t]$.
Let $A = (a_{ij}(t))$ be the full $n\times n$ matrix of coefficients of the 
linear terms, and let $A' = (a_{ij}(t))$, $1 \le i \le n-1$, 
$2 \le j \le n$ be the matrix of coefficients of $x_2,\ldots,x_n$ in the
first $n-1$ equations.  Assume, up to a constant factor,
$\det A'$ is the product of the of the first $n-1$ diagonal entries of the Smith normal form of $A$.
Then solving for $x_2,\ldots,x_n$ from the first $n-1$ equations by 
Cramer's Rule and substituting into the last equation produces
an equation of the form
$$cx_1^2 + B(t)x_1 + C(t) = 0,$$
where $B(t),C(t) \in {\Bbb C}[t]$.  
\endproclaim

\proof
To make the connection between $A$ and $A'$ clearer, we note 
that $A' = A_{n1}$ (submatrix obtained by deleting row $n$ and 
column $1$.  We will number the rows in $A'$ by indices $1$ through $n-1$
and the columns by indices $2$ through $n$ in the following.
As described above for the Selesnick-Burrus equations, take the 
first $n-1$ equations, subtract the $x_1$ terms to the right sides,
and apply Cramer's Rule to solve for $x_2,\ldots, x_n$ in terms
of $x_1$, yielding:
$$x_j = {\det A'_j \over \det A'},$$
for $2 \le j \le n$, where $A'_j$ is the matrix obtained from $A'$ 
by replacing the $j$th column (recall, this means
the column containing the the $a_{ij}(t)$, $ 1 \le i \le n-1$ for this $j$)
with the vector 
$$(r_1(t) - a_{11}(t)x_1,\ldots,r_{n-1}(t)-a_{n-1,1}(t)x_1)^{tr}.\leqno(6.6)$$
If we expand $\det A'_j$ along the column (6.6) in each case we
obtain an expression 
$$\det A'_j = (-1)^{j+1} \det A_{nj} x_1 
+ \sum_{i=1}^{n-1} (-1)^{i+j} r_i(t) \det A'_{ij}$$
where $A_{nj}$ is the submatrix of $A$ obtained by deleting row $n$
and column $j$, and $A'_{ij}$ is a minor of $A'$ (which
is also a submatrix of $A$ obtained by deleting two rows and two columns).
Substitute for $x_i$ in the last equation in the system and rearrange,
taking all the $r_j(t)$ terms to the right hand side.  
The coefficient of $r_j(t)$ is $1/\det A'$ times 
$(-1)^{n+j} \det A_{j1}$.  Hence we obtain
$$
{1\over \det A'}\left(\sum_{j=1}^n a_{nj}\cdot (-1)^{j+1} \det A_{n,j}\right)x_1
= cx_1^2 + {1\over \det A'} \left(\sum_{j=1}^n (-1)^{n+j} \det A_{j1} r_j(t)\right)$$
Up to a sign, the coefficient of $x_1$ is $1/\det A'$ times the
determinant of $A$, expanded along the $n$th row.  
So this can be rewritten as 
$$(-1)^{n-1} {\det A\over \det A'} = c x_1 ^2 + {1\over \det A'} \left(\sum_{j=1}^n (-1)^{n+j} \det A_{j1} r_j(t)\right)$$
By hypothesis, $\det A'$ divides $\det A$ and all of the 
$\det A_{j1}$, which finishes the proof.
\endproof

Now, we must show that the linear equations in the Selesnick-Burrus
system satisfy the hypotheses of the Lemma.  But this follows from 
the determinant formulas from \S 4.
 In our case, $x_1 = m_{2L+3}$, and
$x_2,\ldots, x_{n}$ are $m_{2L+4}, \ldots, m_{3L+4}$.  Continuing 
from Lemma (6.2), $A$ is the matrix $A(L+2,L+1;K)$ (times a lower triangular factor of
determinant 1),
and $A'$ is $\tilde{A}(L+2,L;K)$ (times another lower triangular factor of
determinant 1).  Hence by Theorem (4.10) we have $c_A = 2s + m - 1 + K = 3L+4+K$, and
$$\det A = \delta(L+2,L+1;K) = a \prod_{i=2L+3+K}^{4L+5+K} (2t - i)^{\lfloor {L+1 - |3L+4+K - i|\over 2} \rfloor + 1}$$
for some constant $a$.
Similarly by Corollary (4.12), we have $\tilde{c} = 3L+4 + K$ and 
$$\det A' = \tilde{\delta}(L+2,L;K) = 
a' \prod_{i=2L+5+K}^{4L+3+K} (2t - i)^{\lfloor {L-1 - |3L + 4 + K - i|\over 2} \rfloor + 1},
$$
for some constant $a'$.  Because of the $L - 1$ in the exponent, each factor
in $\det A'$ occurs with multiplicity one less than in $\det A$.
Hence $\det A'$ is precisely the product of the first $L$ diagonal 
entries of the Smith normal form of the $(L+1)\times (L+1)$ matrix $A$, or equivalently, 
${\det A\over \det A'}$ is a polyomial whose roots are all the roots of
$\det A = 0$, but all with multiplicity 1.  Hence the conclusion 
of Lemma (6.5) holds, and we obtain an equation of the form (6.4).

Because of similar cancellations, the final
nonlinear equation (from $F^{(4L+8)}(0) = 0$) has the form
$$[2]m_{2L+3}^2 + {[4L+4]\over [2L-1]} m_{2L+3} + {[6L+7]\over [2L-1]} = 0\leqno(6.6)$$
after we substitute for $m_{2L+4},\ldots,m_{3L+4}$ from (6.3).  The polynomial of 
degree $2L-1$ in the denominators is the same in both terms after the first,
and equals the last diagonal entry in the Smith normal form of 
$\tilde{A}(L+2,L;K)$ (the reduced polynomial of the determinant $\tilde{\delta}(L+2,L;K)$).

The final step is to eliminate $t$ from the two equations (6.4) and (6.6).
For this, we use the determinant form of the Sylvester resultant (see \cite{CLO}) of
two quadratic polynomials, after clearing the denominators in (6.6).  We have
$$\eqalign{Res &= \det\pmatrix [0] & [2L+3] & [4L+6] & 0 \cr
                                0  & [0]    & [2L+3] & [4L+6] \cr
                               [2L+1] & [4L+4]& [6L+7] & 0\cr
                                0  & [2L+1] & [4L+4] & [6L+7] \cr \endpmatrix\cr
               &=  [12L +14]\cr}$$
This concludes the proof of Theorem (6.1).

Computation of these polynomials for a number of $L$ and $K$ suggests
that the polynomial of degree
$[12L + 14]$ is always irreducible over ${\Bbb Q}$, hence has distinct roots in 
${\Bbb C}$.  But we do not have a proof of this fact.  The filters obtained in
this case are considered in \cite{LL}.   

The strategy from (3.10) that we have used here and in \S 4 can also be used to analyze 
the lower diagonals $M = 2L + q$, $q \ge 5$ in Region II.  For instance,
for $q = 5$, solving the linear part of the reduced Selesnick-Burrus system
and substituting into the remaining equations leads to a system of 3 quadrics
in 2 variables (or 3 homogeneous variables).  Explicit determinantal
formulas for the multipolynomial resultant $Res_{2,2,2}$ (see \cite{CLO}, 
Chapter 3, \S 2) can be applied, 
and it can be seen that for $L \ge 2$, the degree of the univariate polynomial
in $t$ is $20L + 26$. We will not present the details of that case
here.

However, the resultants needed 
to eliminate variables in the final, nonlinear system get progressively
harder to analyze as $q$ increases.  Unfortunately, the Dixon resultants 
leading to the most efficient computations tend to have many extraneous factors
that must be accounted for.  As a result, they
are less convenient for the type of analysis done here.   

\refstyle{A}
 
\enddocument